\documentclass[a4paper, 11pt]{amsart}

\usepackage{cite}
\usepackage{graphicx}
\usepackage{color}
\usepackage[T1]{fontenc}

\usepackage{amsmath}
\usepackage{amssymb}
\usepackage{amsthm}
\usepackage{bm}
\usepackage{bbm}
\usepackage{eucal}
\usepackage{braket}
\usepackage{stmaryrd}
\usepackage{mathrsfs}

\usepackage[all,2cell]{xy}
\UseAllTwocells
\usepackage{youngtab}
\usepackage{ytableau}
\ytableausetup{smalltableaux}
\usepackage{tikzit}
\newcommand{\hackcenter}[1]{
 \xy (0,0)*{#1}; \endxy}
 \calclayout

\theoremstyle{plain}
\newtheorem{thm}{Theorem}[section]
\newtheorem{lem}[thm]{Lemma}
\newtheorem{prop}[thm]{Proposition}
\newtheorem{cor}[thm]{Corollary}

\theoremstyle{definition}
\newtheorem{defn}[thm]{Definition}%[section]

\theoremstyle{remark}
\newtheorem{rem}[thm]{Remark}%[section]
\newtheorem{exam}[thm]{Example}
\makeatletter
    
    \@addtoreset{equation}{section}
\makeatother

\usepackage{xr-hyper}
\usepackage{hyperref} %comment out in setup.tex

\makeatletter

\newcommand*{\addFileDependency}[1]{% argument=file name and extension
\typeout{(#1)}% latexmk will find this if $recorder=0
\@addtofilelist{#1}
\IfFileExists{#1}{}{\typeout{No file #1.}}
}\makeatother

\newcommand{\bsl}{\backslash}

\newcommand{\wtilde}[1]{\widetilde{#1}}
\newcommand{\ul}[1]{\underline{#1}}

\newcommand{\End}{\mathrm{End}}

\newcommand{\Tr}{{\mathrm{{Tr}}}}

\newcommand{\id}{\mathrm{id}}

\newcommand{\catMod}{\mathsf{Mod}}

\newcommand{\funInd}{\mathsf{Ind}}
\newcommand{\funRes}{\mathsf{Res}}
\newcommand{\funId}{\mathsf{Id}}

\newcommand{\bC}{\mathbb{C}}

\newcommand{\bM}{\mathbb{M}}
\newcommand{\bN}{\mathbb{N}}

\newcommand{\bR}{\mathbb{R}}

\newcommand{\bY}{\mathbb{Y}}
\newcommand{\bZ}{\mathbb{Z}}

\newcommand{\cD}{\mathcal{D}}
\newcommand{\cE}{\mathcal{E}}
\newcommand{\cF}{\mathcal{F}}

\newcommand{\cH}{\mathcal{H}}
\newcommand{\cI}{\mathcal{I}}

\newcommand{\cP}{\mathcal{P}}

\newcommand{\sfm}{\mathsf{m}}
\newcommand{\x}{\mathsf{x}}
\newcommand{\y}{\mathsf{y}}

\newcommand{\bbmi}{\mathbbm{i}}

\newcommand{\bfE}{\mathbf{E}}
\newcommand{\bfT}{\mathbf{T}}

\makeatletter
\newcommand{\xRightarrow}[2][]{\ext@arrow 0359\Rightarrowfill@{#1}{#2}}
\makeatother

\newcommand{\ST}{\mathrm{ST}}
\newcommand{\Conf}{\mathrm{Conf}}
\newcommand{\Pl}{\mathrm{Pl}}

\newcommand{\wt}{\mathrm{wt}}

\newcommand{\ex}{\mathrm{ex}}
\tikzstyle{dot}=[fill=black, draw=black, shape=circle, minimum size=10pt, inner sep=0pt]
\tikzstyle{ExtPt}=[fill={rgb,255: red,128; green,128; blue,128}, draw={rgb,255: red,128; green,128; blue,128}, shape=circle, minimum size=5pt, inner sep=0pt]
\tikzstyle{Sdot}=[fill=black, draw=black, shape=circle, minimum size=5pt, inner sep=0pt]

\tikzstyle{thick arrow}=[->, thick]
\tikzstyle{thick edge}=[-, thick]
\tikzstyle{verythick edge}=[-, very thick]
\tikzstyle{normal edge}=[-]
\tikzstyle{dashed line}=[-, dashed]
\tikzstyle{UTRED}=[-, ultra thick, draw=red]
\tikzstyle{shade}=[-, fill={rgb,255: red,100; green,255; blue,255}, draw={rgb,255: red,128; green,128; blue,128}, ultra thick]
\tikzstyle{dotted arrow}=[->, dotted, thick]
\tikzstyle{verythickArrow}=[very thick, ->]
\tikzstyle{ExtDisk}=[-, ultra thick, draw={rgb,255: red,128; green,128; blue,128}]
\tikzstyle{GRShade}=[-, fill={rgb,255: red,128; green,128; blue,128}, draw={rgb,255: red,128; green,128; blue,128}]

\title[Planar algebra and Heisenberg category]{Planar algebras for the Young graph and the Khovanov Heisenberg category}
\date{\today}
\author{Shinji Koshida} %% Under amsart
\address{Department of Mathematics and Systems Analysis, Aalto University, Finland} %% Under amsart
\email{shinji.koshida@aalto.fi} %% Under amsart

\begin{document}

\maketitle

\begin{abstract}
This paper studies planar algebras of Jones' style associated with the Young graph.
We first see that, given a positive real valued function on the Young graph,
we may obtain a planar algebra whose structure is defined in terms of a state sum
over the ways of filling planar tangles with Young diagrams.
We delve into the case that the function is harmonic and related to the Plancherel measures on Young diagrams.
Along with an element that is depicted as a cross of two strings,
we see that the defining relations among morphisms for the Khovanov Heisenberg category
are recovered in the planar algebra.
We also identify certain elements in the planar algebra with particular functions of Young diagrams
that include the moments, Boolean cumulants and normalized characters.
This paper thereby bridges diagramatical categorification and asymptotic representation theory.
In fact, the Khovanov Heisenberg category is one of the most fundamental examples of diagramatical categorification
whereas the harmonic functions on the Young graph have been a central object in the asymptotic representation theory
of symmetric groups.
\end{abstract}

\tableofcontents

\section{Introduction}

\subsection{Rattan--{\'S}niady conjecture}
This paper can be considered a companion paper of the previous one of the author's~\cite{Koshida2023}.
In order to see how the present work came up, let us start off by the Rattan--{\'S}niady conjecture.

The conjecture is concerned with two series of functions of Young diagrams; normalized characters originating from the representation theory of symmetric groups, and Boolean cumulants that measure the shape of Young diagrams.
We will give precise definitions of those functions in Section~\ref{Sect:young}, and only put few words here.

They are both functions on the set of Young diagrams $\bY$.
The normalized characters denoted by $\Sigma_{\pi}$ are labelled by partitions $\pi\in\cP$. The set of Young diagrams $\bY$ and that of partitions $\cP$ are naturally identified, but it is often mentally useful to distinguish them according to their different roles.
The Boolean cumulants $B_{k}$ are, in principle, labelled by positive integers $k\in\bZ_{> 0}$, but in the context of Young diagrams, we always have $B_{1}=0$. Thus, $B_{k}$, $k\geq 2$ are the only non-trivial functions.

It is not even clear from the definition that each normalized character is expressed as a polynomial of Booelan cumulants, but it was shown to be the case in~\cite{Biane2003}.
The conjecture by Rattan and {\'S}niady~\cite{RS2008}, now a theorem due to the present author, states more:
\begin{thm}[\cite{Koshida2023}, conjectured in~\cite{RS2008}]
\label{thm:RSconj}
For each $\pi\in \cP$, there is a polynomial $P_{\pi}(x_{2},\dots,x_{|\pi|-\ell(\pi)+2})$ with $\bZ_{\geq 0}$-coefficients such that
\begin{align*}
    (-1)^{\ell(\pi)}\Sigma_{\pi} = P_{\pi}(-B_{2},\dots,-B_{|\pi|-\ell(\pi)+2}).
\end{align*}
Here, $|\pi|$ and $\ell (\pi)$ are the weight and length of the partition $\pi$.
\end{thm}

Our method of proving Theorem~\ref{thm:RSconj} was, however, very indirect;
we used diagram calculus in the Khovanov Heisenberg category~\cite{Khovanov2014}.

\subsection{Khovanov Heisenberg category}
\begin{figure}
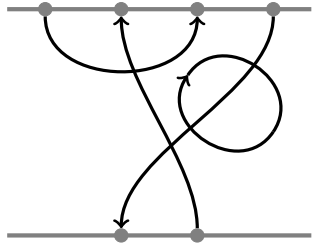

\ctikzfig{khMorphism}
\caption{A morphism from $Q_{-}\otimes Q_{+}$ to $Q_{-}\otimes Q_{+}\otimes Q_{+}\otimes Q_{-}$.}
\label{fig:HeisenbergCatMorph}
\end{figure}

The Khovanov Heisenberg category, denoted by $\cH$, is a strict monoidal category freely generated by two objects $Q_{+}$ and $Q_{-}$.
The morphisms from $Q_{\epsilon_{1}}\otimes\cdots \otimes Q_{\epsilon_{k}}$ to $Q_{\epsilon'_{1}}\otimes \cdots \otimes Q_{\epsilon'_{l}}$ are described as follows.
We place points with the signs $\epsilon_{1},\dots, \epsilon_{k}$ on $\bR\times \{0\}\subset \bR^{2}$ in this order and
points with the signs $\epsilon'_{1},\dots,\epsilon'_{l}$ on $\bR\times\{1\}\subset\bR^{2}$ similarly.
A morphism is a collection of oriented strings in $\bR\times [0,1]$, each of which either connects two of the points respecting their signs or forms a closed loop (Figure~\ref{fig:HeisenbergCatMorph}). We require that each point is connected to exactly one string, and that the strings have only normal crossings.
We also consider the strings only up to isotopy, but most importantly, they are subject to the following local relations:
\begin{align}
\label{eq:KHrelation1}
	\input{khIndIndCrossings.tikz}\quad  = \quad \input{khIndIndStraight.tikz}\quad ,\qquad
	\input{khIndResCrossings.tikz} \quad = \quad \input{khIndResStraight.tikz}\quad ,
\end{align}
\begin{align}
\label{eq:KHrelation2}
	\input{khResIndCrossings.tikz}\quad = \quad \input{khResIndStraight.tikz}\quad -\quad \input{khResIndCups.tikz}\quad ,
	\qquad \input{khYBEleft.tikz}\quad = \quad \input{khYBEright.tikz}\quad ,
\end{align}
\begin{align}
\label{eq:KHrelation3}
	\input{khLeftTurn.tikz}\quad = 0,\qquad \input{khLeftCircle.tikz}\quad =1.
\end{align}

\begin{rem}
In most cases in the context of categorification, the Khovanov Heisenberg category is the Karoubi envelop of the above category,
but for us the above definition suffices.
\end{rem}

\subsection{Proof of Rattan--{\'S}niady conjecture}
The strategy of proving Theorem~\ref{thm:RSconj} in~\cite{Koshida2023} was as follows.
First of all, the Khovanov Heisenberg category has the following representation~\cite{Khovanov2014} (and this is probably how the category was discovered)
\begin{align*}
	\cF \colon \cH \to \cE nd\Big( \bigoplus_{n\geq 0}\bC [S_{n}]-\catMod \Big)
\end{align*}
such that
\begin{align*}
	\cF (Q_{+}) = \bigoplus_{n\geq 0} \funInd_{S_{n}}^{S_{n+1}},\quad \cF (Q_{-})=\bigoplus_{n\geq 0} \funRes^{S_{n+1}}_{S_{n}}.
\end{align*}
This functor, in particular, induces an algebra homomorphism
\begin{align*}
	\End_{\cH}(\bm{1}) \to \End_{\cE nd \big( \bigoplus_{n\geq 0}\bC [S_{n}]-\catMod\big)}(\funId),
\end{align*}
where $\bm{1}$ is the unit object of $\cH$.
Now, if we expand the definition of natural transformations on the identity functor,
we can identify the algebra in the right-hand side with $\prod_{n\geq 0}Z(\bC [S_{n}])$,
which, under the Fourier transform, can be identified further with $\bC [\bY]$, the algebra of functions of Young diagrams.

Composing all these identifications, we get an algebra homomorphism
\begin{align}
\label{eq:KLMhom}
	\End_{\cH}(\bm{1}) \to \bC [\bY]
\end{align}
out of the functor $\cF$.
The dictionary due to~\cite{KLM2019} gives us the diagrams in the left-hand side that are sent to the normalized characters and the Boolean cumulants.
Thus, we can study algebraic relation among those functions in terms of diagram calculus in the left-hand side, which is exactly what we did in the previous work~\cite{Koshida2023} to prove Theorem~\ref{thm:RSconj}.

\subsection{Questions}
When we observe non-negative integers in the coefficients of an expansion as in Theorem~\ref{thm:RSconj},
what is often the case is that those integers are counting some combinatorial objects.
Thus, it is ideal to have a combinatorial interpretation of the coefficients in our expansion as well.

So far, our proof of  Theorem~\ref{thm:RSconj} seems hopeless in this regard;
it is a combination of diagram calculus and induction, but the induction procedure we used is too complicated to see a combinatorial interpretation of the expansion coefficients.
One of the reasons for the complication is probably that the induction involves diagrams that are neither normalized characters nor Boolean cumulants.

Nonetheless, we can still ask what combinatorics could possibly happen behind the story.
In particular, we seek a combinatorial understanding of our main tools: the homomorphism (\ref{eq:KLMhom}) and the dictionary of~\cite{KLM2019}.
In turns out that the question fits the setting of planar algebras of Jones' style~\cite{jones2021planar}, and from there, we found a way to rediscover the Khovanov Heisenberg category.

\subsection{Planar algebras {\`a} la Jones}
Planar algebras were originally introduced in the context of subfactors.
We do not give a precise definition, but roughly, a planar algebra consists of a collection of vector spaces along with a collection of multi-linear maps that are labelled by the so-called planar tangles.
The axioms for a planar algebra require those collections to be invariant under isotopy and consistent under composition (see~\cite{jones2021planar} for details).

A prominent construction of a planar algebra starts from a bipartite graph~\cite{jones2000planar}. It actually depends not only on a bipartite graph, but also on a non-negative eigenfunction of the adjacency matrix of the graph.
Then, we collect the spaces of functions of loops in the graph of various lengths, and define multi-linear maps in terms of sums over possible fillings of the regions of planar tangles with vertices of the graph.

\subsection{Planar algebras for the Young graph}
The idea of {\it filling the regions with vertices} turns out helpful for our purpose of understanding the Khovanov Heisenberg category combinatorially.
In fact, we are going to plug the Young graph (see Section~\ref{Sect:young}) into the place of a bipartite graph and get planar algebras.
However, our case does not really fall into a special case of Jones' planar algebras.

In the case of a bipartite graph, we need a positive eigenfunction of the adjacency matrix of the graph,
whose existence is guaranteed by the Perron--Frobenius theorem.
On the other hand, it is natural to think of the Young graph from its origin as an {\it oriented} graph.
Thus, its adjacency matrix is no longer symmetric.
Luckily for us, there is a classification of the harmonic functions on the Young graph~\cite{VershikKerov1981, Kerov2003, BorodinOlshanski_book2017},
so it seems to be a good option to take a harmonic function, i.e., a positive eigenfunction of the adjacency matrix of eigenvalue $1$,
as an input for our planar algebra.

The fact that the Young graph is oriented gives rise to another feature of our planar algebra, that is,
the planar tangles that label multi-linear maps should also be oriented.
This is, of course, natural since the strings in the Khovanov Heisenberg category are oriented too.

Strictly speaking, Jones' construction of a planar algebra out of a bipartite graph~\cite{jones2000planar} does not require the function to be an eigenfunction of the adjacency matrix,
but assuming that it is an eigenfunction allows for a further property, namely, loop removing.
We will see in Section~\ref{Sect:planarAlg} that, given a positive real valued function on the Young graph, we get a planar algebra.
The vector spaces for it are the spaces of functions of loops in the Young graph (forgetting the orientation) of various lengths
and the multi-linear maps are defined by summing over fillings of the regions of planar tangles with Young diagrams.
We also observe loop removing under the assumption of harmonicity of the input function.

\subsection{Plancherel case}
While there are a lot of harmonic functions on the Young graph,
there is a particularly prominent one; the harmonic function coming from the Plancherel measures.

We will relate the planar algebra associated with this specific choice of harmonic function to the Khovanov Heisenberg category,
but for that, we need one further ingredient, namely, crossing.
By definition of a planar tangle, it contains no crossing of strings whereas in the Khovanov Heisenberg category, strings can cross. 
We introduce a crossing as a special element of a vector space that is part of the planar algebra.
More specifically, for us, a crossing is a function of loops in the Young graph of length four.

Under this setup, we will see that the relations (\ref{eq:KHrelation1})--(\ref{eq:KHrelation3}) imposed on the morphisms in the Khovanov Heisenberg category are recovered as relations in the planar algebra.
Furthermore, various functions of Young diagrams including the normalized characters and the Boolean cumulants are realized as elements of the planar algebra.
As such, we could say that the Khovanov Heisenberg category and the dictionary of \cite{KLM2019} are rediscovered from the planar algebra.

\subsection{Discussions}
\subsubsection{Interplay between planar algebras, combinatorics and others}
As we have already mentioned, this paper can be considered a companion paper of our previous one~\cite{Koshida2023}.
This is in the following sense. In the previous work, we used the homomorphism (\ref{eq:KLMhom}) and diagram calculus in $\End_{\cH}(\bm{1})$ to prove Theorem~\ref{thm:RSconj}.
The present work allows us to perform the diagram calculus already in $\bC[\bY]$; the space of functions on $\bY$ is the same as
the space of functions of loops in $\bY$ of length $0$, which is part of our planar algebra.
Thus, we could prove Theorem~\ref{thm:RSconj} starting from the planar algebra, but without the Khovanov Heisenberg category.

On the other hand, we could reverse this point to ask another question.
Just for proving Theorem~\ref{thm:RSconj}, we only need either our planar algebra, or the version of the Khovanov Heisenberg category before the Karoubi enveloping.
In the context of categorification, however, it seems that a category is well-behaved if it is Karoubi complete.
One could then ask if the usual, Karoubi completed, Heisenberg category can be used to answer some questions
from asymptotic representation or combinatorics of Young diagrams.

Another way to see the present work, is that it is bridging between
diagramatial categorification and asymptotic representation by language of planar algebras.
In fact, on the one hand, the Khovanov Heisenberg category provides one of the most prominent categorifications of algebraic structures, and on the other hand, the Young graph, in particular harmonic analysis on it, is in the center of
the asymptotic representation of symmetric groups.
Khovanov~\cite{Khovanov2014} already pointed out
\begin{quote}
``It also seems that our construction should be related to the circle of ideas considered by Guionnet, Jones and Shlyakhtenko~\cite{GuionnetJonesShlyakhtenko2010} that intertwine planar algebras and free probability.''
\end{quote}
(Change of citation index by the author.)
The present work makes a part of the comment clearer giving a precise definition of the planar algebra that pertains to the Khovanov Heisenberg category.
We did not discuss free probability, though, and it is a fair question for the future what a random matrix/free probability model computes the trace of our planar algebra.

\subsubsection{General branching graphs}
The only reason of taking the Young graph as input in this paper is because the outcome recovers the Khovanov Heisenberg category.
From a more general perspective of asymptotic representation theory,
we should, of course, be able to start from a general branching graph (see~\cite{petrov2013operators} for possibly the most general notion of branching graph).
The Khovanov Heisenberg category also admits variants and deformations~\cite{LicataSavage2013, MackaaySavage2018,Brundan2018, BrundanSavageWebster2020}.
Extending the correspondence appearing in this paper to other cases will be an important direction to establishing the relationship between diagram categories and asymptotic representations.

\subsubsection{Jack multiplicities}
But perhaps, there is more than just the asymptotic representation of some algebraic structure.
As long as our input data are a branching graph and a harmonic function, there is another direction of generalization:
introducing weights on the edges.
It is common to introduce weights to the Young graph that originate from the Pieri rules of symmetric functions,
and harmonic analysis with weights works well, particularly with the Jack weights~\cite{kerov1998boundary, kerov2000anisotropic} (and the Macdonald weights~\cite{matveev2019macdonald}).
Analyzing the Jack characters is, however, more difficult than the symmetric group case because of the lack an actual representation theory underlying.
It seems likely that our construction of planar algebras admits Jack deformation if not entirely straightforward,
and if it does, it will provide a new approach to the Jack deformed {\it asymptotic representation}~\cite{Lassalle2009,DolegaFeraySniady2014,dolkega2016gaussian,sniady2019asymptotics},
and possibly a diagramatical categorification of Jack polynomials.
Let us mention~\cite{LicataRossoSavage2018} that categorified the Jack inner product using the Frobenius Heisenberg categorification~\cite{CautisLicata2012, RossoSavage2017, Savage2019}, but our possible approach from the planar algebra seems different.

\subsubsection{Beyond the Plancherel case}
There is also a remaining question in the Young graph case.
Although we have studied the planar algebra for the Young graph with the Plancherel harmonic function,
we know less about other choices of harmonic functions.
For example, there are harmonic functions that are related to $z$-measures, and these functions are well-studied in asymptotic representation theory~\cite{kerov2004harmonic}.
We could then ask if the associated planar algebras are used to discover new diagramatic categories or new representations of the Khovanov Heisenberg category.

\subsubsection{Free cumulants}
In the context of asymptotic representation, free cumulants of Young diagrams often behave better than the Boolean cumulants~\cite{Biane1998}. Also the analogue of Theorem~\ref{thm:RSconj} for the free cumulants is more prominent as the Kerov conjecture, and has been proved in~\cite{Feray2009, DolegaFeraySniady2010} with an explicit combinatorics.
It is unclear, however, how the free cumulants are realized in our planar algebra.
Perhaps, we need to take into account more nonlocal effects than just a crossing of two strings.
Let us remark that there are analogues of the Kerov conjecture under the Jack-deformed~\cite{Lassalle2009} and spin settings~\cite{Matsumoto2018,MatsumotoSniady2020}.

\subsubsection{Contribution of the present work}
If we focus on the Plancherel harmonic function, our results can be, in principle, derived by
carefully expanding the contents of the homomorphism (\ref{eq:KLMhom}) given that we know how the resulting planar algebra looks like.
However, it is arguably nontrivial to find the correct setting for a planar algebra. Also, we would imagine that deriving the relations (\ref{eq:KHrelation1})--(\ref{eq:KHrelation3}) in our approach is not so common in the community of categorification.
On the other hand, applying diagram calculus through the homomorpshim (\ref{eq:KLMhom}) that is induced from a functor is alien in the community of asymptotic representation.
For those reasons, we believe that it is still worthwhile recording the contribution of the present paper.

Finally, let us comment that the original problem of finding combinatorics behind Theorem~\ref{thm:RSconj} is still open.

\subsection*{Organization}
This paper is organized as follows.
In Section~\ref{Sect:young}, we recall known facts about Young diagrams and the Young graph that are needed in the rest of the paper.
In Section~\ref{Sect:planarAlg}, we fix the definition of a planar algebra.
Although the notion therein is standard, we make a complete presentation as our version of a planar algebra might differ from other literatures in small details.
We also construct a planar algebra associated to the Young graph and a positive real valued function on it.
In Section~\ref{Sect:plancherelCase}, we study the specific example when we take the harmonic function that originates from the Plancherel measures.
We observe the defining relations (\ref{eq:KHrelation1})--(\ref{eq:KHrelation3}) of the Khovanov Heisenberg category and the dictionary due to~\cite{KLM2019} in our planar algebra.
This paper contains two appendices.
Appendix~\ref{app:reptheorSG} recollects some facts about representations of symmetric groups
that are needed to complete the proof of Theorem~\ref{thm:characterDiagram} in Section~\ref{Sect:plancherelCase}.
In Appendix~\ref{app:charactersAlter}, we seek an alternative way of proving Theorem~\ref{thm:characterDiagram}
that does not require representation theory.

\subsection*{Acknowledgments}
The author is grateful to Piotr {\'S}niady for the discussions that gave rise to the initial idea for the present work.
The author also thanks the anonymous referee for pointing out errors present in the first submission and various suggestions for improvement.
This work was supported by funding from Academy of Finland (No. 340965).

\section{Young diagrams, graph and characters}
\label{Sect:young}
This section collects known facts regarding Young diagrams and the Young graph
and introduce various functions of Young diagrams including the moments, Boolean cumulants and normalized characters.
More details can be found in~\cite{Fulton1996, kerov2000anisotropic, Biane2003}.

\subsection{Young diagrams}
For us, a Young diagram is a non-increasing sequence of non-negative integers $\lambda = (\lambda_{1},\lambda_{2},\dots)$
such that the sum $|\lambda| = \sum_{i\geq 1}\lambda_{i}$, called the weight, is finite.
In other words, $\lambda_{i}$'s are all $0$ except for the first few.
We often write only non-zero entries of a Young diagram as $\lambda = (\lambda_{1},\dots,\lambda_{l})$.
%The number of non-zero entires is denoted by $\ell (\lambda) = l$ and called the length of the partition $\lambda$.
 %\cyan{Check if this terminology is really needed.}

A Young diagram is a {\it diagram} because it is usually drawn as a pile of as many boxes as its weight.
There are a few conventions of doing this, but one way is such that the number of boxes is non-increasing downward (English convention).
For example, we can draw
\begin{equation*}
    (7,4,2,1) = 
    \hackcenter{\begin{ytableau}
        \, & \, & \, & \, & \, &\, &\, \\
        \, & \, & \, & \, \\
        \, & \, \\
        \,
    \end{ytableau}}.
\end{equation*}
In the following, we do not distinguish a Young diagram as a sequence of non-negative integers from
its graphical presentation if not necessary.

When we reflect a Young diagram drawn in the above convention with respect to
the straight line with slope $-1$ passing through the left-top corner,
we get another Young diagram.
We call this Young diagram the transpose of the original diagram.
For a Young diagram $\lambda$, its transpose will be denoted by $\lambda'$.
For example, the transpose of $\lambda = (7,4,2,1)$ (the above example) is
\begin{equation*}
	\lambda' = (4,3,2,2,1,1,1) =
	\hackcenter{\begin{ytableau}
        \, & \, & \, & \, \\
        \, & \, & \,  \\
        \, & \, \\
        \, & \, \\
        \, \\
        \, \\
        \,
    \end{ytableau}}\quad .
\end{equation*}

Each box of a Young diagram is identified in the format $(i,j)$ of a pair of positive integers,
where the first and second components indicate the downward and rightward coordinates, respectively.
For example, the shaded box in
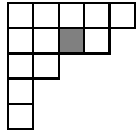
\begin{equation*}
    \hackcenter{\begin{ytableau}
        \, & \, & \, & \, & \, & \, & \, \\
        \, & \, & *(gray) & \, \\
        \, & \, \\
        \,
    \end{ytableau}}
\end{equation*}
has the coordinates $(2,3)$.
It is often convenient to write $(i,j)\in \lambda$ when the box identified by $(i,j)$ is in the Young diagram $\lambda$.
For a box $(i,j)$, its content is defined by
\begin{equation*}
	c((i,j)) := j - i.
\end{equation*}
In other words, it indicates the distance from the diagonal.

For two Young diagrams $\mu$ and $\lambda$, we say that $\mu$ is contained in $\lambda$ and write $\mu\subset\lambda$
if $\mu_{i}\leq\lambda_{i}$ for all $i$.
In this case, the difference of the two Young diagrams forms a skew diagram denoted by $\lambda/\mu$.

\subsection{Young graph}
\label{subsect:YoungGraph}
\begin{figure}
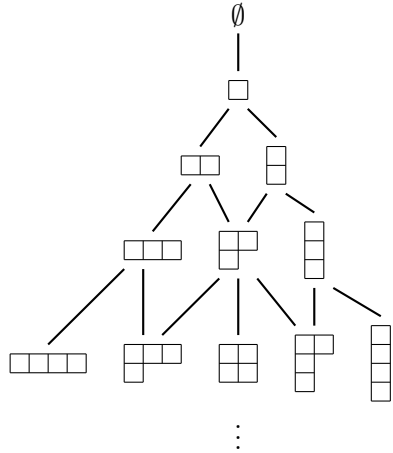

\ctikzfig{youngGraph}
\caption{Young graph. The vertical direction (downwards) indicates the weights of Young diagrams.}
\label{fig:YoungGraph}
\end{figure}

For $n\in \bN$, the set of Young diagrams of weight $n$ will be denoted by $\bY_{n}$.
We set $\bY_{0}=\{\emptyset\}$ by convention.
For $\lambda\in \bY_{n+1}$ and $\mu\in \bY_{n}$ with $n\in\bN$, we write $\mu\nearrow \lambda$ if $\mu\subset \lambda$, i.e., $\lambda$ is obtained by adding one box to $\mu$.
As such, we get the Young graph with the vertex set $\bY:=\bigcup_{n\in\bN}\bY_{n}$ and edges $\nearrow$ (Figure~\ref{fig:YoungGraph}).
When $\mu\nearrow \lambda$, we might also write $\lambda\searrow\mu$ depending on convenience.

When $\mu\nearrow\lambda$, the skew diagram $\lambda/\mu$ consists of a single box.
It is reasonable, then, to write $c(\lambda/\mu)$ for the content of the box.

A path in the Young graph is a sequence $(\lambda^{(0)},\lambda^{(1)},\dots,\lambda^{(n)})$ of Young diagrams such that
\begin{align*}
	\lambda^{(0)}\nearrow \lambda^{(1)}\nearrow \cdots \nearrow \lambda^{(n)}.
\end{align*}
When we talk about a loop, we need to go back along edges,
but it is often convenient to remember in which direction we go through each edge.
For a sign $\epsilon\in\{+,-\}$, we introduce the following notation
\begin{align*}
	\lambda \xrightarrow{\epsilon} \mu :=
	\begin{cases}
		\lambda \nearrow \mu, & \epsilon = +, \\
		\lambda \searrow \mu, & \epsilon = -,
	\end{cases}
\end{align*}
for two Young diagrams $\lambda$ and $\mu$.
Let $n\in \bN$ and $\ul{\epsilon} = (\epsilon_{1},\dots,\epsilon_{n})$ be a sequence of signs, i.e., $\epsilon_{i}\in \{+,-\}$, $i=1,\dots, n$ such that $\sum_{i=1}^{n}\epsilon_{i}=0$.
A loop in the Young graph of signature $\ul{\epsilon}$ is a sequence $(\lambda^{(0)},\lambda^{(1)},\dots,\lambda^{(n)})$ of Young diagrams such that
\begin{align*}
	\lambda^{(0)}\xrightarrow{\epsilon_{1}} \lambda^{(1)} \xrightarrow{\epsilon_{2}} \cdots \xrightarrow{\epsilon_{n}} \lambda^{(n)} = \lambda^{(0)}.
\end{align*}
In both cases of a path and a loop, we call the number $n$ appearing above the length.

The following observation for paths of length $2$ will be useful in many places.
Let $\lambda \nearrow \mu \nearrow \nu$ be a path in the Young graph.
There are two qualitatively distinguishable cases.
The first case is that $c(\nu/\mu) - c(\mu/\lambda) = \pm 1$.
In this case the two boxes $\mu/\lambda$ and $\nu/\mu$ are added in the same row or column.
An example of such a case is
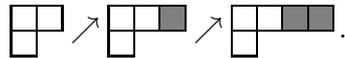
\begin{align*}\,
	\hackcenter{\begin{ytableau}
        \, & \, \\
        \, 
        \end{ytableau}}
        \nearrow 
        \hackcenter{\begin{ytableau}
        \, & \, & *(gray) \\
        \, 
        \end{ytableau}}
        \nearrow\hackcenter{\begin{ytableau}
        \, & \, & *(gray) & *(gray) \\
        \, 
        \end{ytableau}}\, .
\end{align*}
The other case is that $c(\nu/\mu) - c(\mu/\lambda) \neq \pm 1$.
In this case, we can change the order of adding the two boxes $\mu/\lambda$ and $\nu/\mu$
and get another path $\lambda\nearrow \mu'\nearrow \nu$.
For example, we have the following pair of paths:
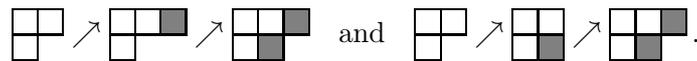
\begin{align*}\,
	\hackcenter{\begin{ytableau}
        \, & \, \\
        \, 
        \end{ytableau}}
        \nearrow 
        \hackcenter{\begin{ytableau}
        \, & \, & *(gray) \\
        \, 
        \end{ytableau}}
        \nearrow\hackcenter{\begin{ytableau}
        \, & \, & *(gray) \\
        \, & *(gray)
        \end{ytableau}}\quad \text{and}\quad
        \hackcenter{\begin{ytableau}
        \, & \, \\
        \, 
        \end{ytableau}}
        \nearrow 
        \hackcenter{\begin{ytableau}
        \, & \,  \\
        \, & *(gray)
        \end{ytableau}}
        \nearrow\hackcenter{\begin{ytableau}
        \, & \, & *(gray) \\
        \, & *(gray)
        \end{ytableau}}\, .
\end{align*}

%\cyan{Delete signature explanation for loops later. Convention for $\xrightarrow{\epsilon}$ has also changed.}

\subsection{(Co)transition measures}
\begin{figure}
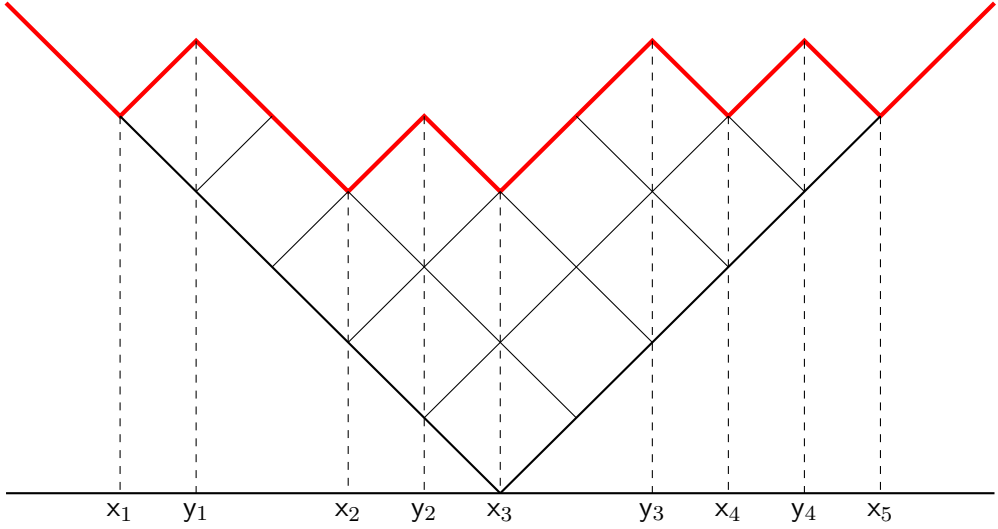

\ctikzfig{youngDiagramRus}
\caption{Profile (red) of the Young diagram $\lambda = (7,4,2,1)$. In this example, the local minima and maxima are $(\x_{1},\x_{2},\x_{3},\x_{4},\x_{5})=(-4,-2,0,3,7)$ and $(\y_{1},\y_{2},\y_{3},\y_{4})= (-3,-1,2,6)$.}
\label{fig:profile}
\end{figure}

There is another convention to draw a Young diagram,
in which we reflect a Young diagram drawn in the above convention vertically and rotate it by 45$^{\circ}$ (Russian convention).
We place the corner at the origin and set the linear length of each box to be $\sqrt{2}$.
Then, we find the profile of the Young diagram that envelopes the boxes along with the lines of slope $\pm 1$ (Figure~\ref{fig:profile}).
We would suggest~\cite{kerov1997interlacing} for a canonical reference to this picture of Young diagrams.

Let us fix a Young diagram $\lambda$.
In the profile of $\lambda$, we get the local minima $(\x_{1},\x_{2},\dots,\x_{d})$ and the local maxima $(\y_{1},\y_{2},\dots, \y_{d-1})$, which are interlacing
\begin{equation}
\label{eq:interlacing}
	\x_{1}<\y_{1}<\x_{2}<\y_{2}<\cdots <\y_{d-1}<\x_{d}
\end{equation}
and the sums of $\x_{i}$'s and $\y_{i}$'s coincide:
\begin{equation}
\label{eq:maxmin_centered}
	\sum_{i=1}^{d}\x_{i} = \sum_{i=1}^{d-1}\y_{i}.
\end{equation}
As we set the linear length of a box to be $\sqrt{2}$, the local maxima and minima are all integers.
Note that the $\x_{i}$'s are the contents of those boxes that are addable to $\lambda$,
and the $\y_{i}$'s are those of removable boxes from $\lambda$.

We then form the rational function
\begin{equation}
\label{eq:def_of_G_rational}
	G_{\lambda}(z) = \frac{\prod_{i=1}^{d-1}(z-\y_{i})}{\prod_{i=1}^{d}(z-\x_{i})},
\end{equation}
which is called the (moment) generating function for $\lambda$~\cite{kerov1993transition,Kerov2003}.
Because of the interlacing structure (\ref{eq:interlacing}), the function $G_{\lambda}$ is the Cauchy transform of a probability measure $\sfm_{\lambda}$ on $\bR$:
\begin{equation*}
	G_{\lambda}(z) = \int_{\bR}\frac{1}{z-x}\sfm_{\lambda}(dx).
\end{equation*}
More explicitly, we can decompose $G_{\lambda}(z)$ into partial fractions as
\begin{equation}
\label{eq:Cauchy_partial_fractions}
	G_{\lambda}(z) = \sum_{i=1}^{d}\frac{1}{z-\x_{i}}\frac{\prod_{j=1}^{d-1}(\x_{i}-\y_{j})}{\prod_{j=1;j\neq i}^{d}(\x_{i}-\x_{j})},
\end{equation}
from which we can see that the probability measure $\sfm_{\lambda}$ has the form
\begin{equation*}
	\sfm_{\lambda} = \sum_{i=1}^{d}\frac{\prod_{j=1}^{d-1}(\x_{i}-\y_{j})}{\prod_{j=1;j\neq i}^{d}(\x_{i}-\x_{j})} \delta_{\x_{i}}.
\end{equation*}

%Let us note that $x_{i}$'s are the contents of the boxes that are addable to $\lambda$
%to get another Young diagram with one more box.
Suppose that another Young diagram $\mu$ is such that $\lambda\nearrow \mu$ and $\x_{i} = c(\mu/\lambda)$ with some $i$.
Then, we set
\begin{equation}
\label{eq:transition_rational}
	p^{\uparrow}(\lambda,\mu) = \frac{\prod_{j=1}^{d-1}(\x_{i}-\y_{j})}{\prod_{j=1;j\neq i}^{d}(\x_{i}-\x_{j})},
\end{equation}
which is the weight put on $\delta_{\x_{i}}$ in $\sfm_{\lambda}$.
We stress that $(\x_{1},\dots, \x_{d})$ and $(\y_{1},\dots, \y_{d-1})$ in this formula are
the local maxima and minima of the profile of $\lambda$, but {\bf not} that of $\mu$.
We have
\begin{equation*}
	\sum_{\mu;\lambda\nearrow\mu}p^{\uparrow}(\lambda,\mu)=1
\end{equation*}
because $\sfm_{\lambda}$ is a probability measure.
This means that $p^{\uparrow}$ gives a stochastic transition on the Young graph $\bY$.
Because of this reason, $\sfm_{\lambda}$ is called the transition measure of $\lambda$~\cite{kerov1993transition}.
Depending on the context, we also call the probability measure on $\bY_{|\lambda|+1}$
defined by $p^{\uparrow}(\lambda,-)$ the transition measure of $\lambda$.
The transition measure was used to study the asymptotic representation of the symmetric groups.
For detailed account and its application, see~\cite{Kerov2003}

It will be convenient to write (\ref{eq:Cauchy_partial_fractions}) in the following form:
\begin{equation}
\label{eq:Cauchy_contents}
	G_{\lambda}(z) = \sum_{\mu;\lambda\nearrow\mu}\frac{p^{\uparrow}(\lambda,\mu)}{z-c(\mu/\lambda)}.
\end{equation}

We can find another probability measure associated with a Young diagram.
We are still working with the fixed Young diagram $\lambda$.
Let us take the multiplicative inverse of $G_{\lambda}(z)$:
\begin{equation}
\label{eq:CauchyInv_rational}
	H_{\lambda}(z):=\frac{1}{G_{\lambda}(z)} = \frac{\prod_{i=1}^{d}(z-\x_{i})}{\prod_{i=1}^{d-1}(z-\y_{i})}.
\end{equation}
Here, $(\x_{1},\dots,\x_{d})$ and $(\y_{1},\dots, \y_{d-1})$ are again the local maxima and minima of the profile of $\lambda$.
This time $H_{\lambda}$ is not the Cauchy transform of a probability measure,
but after extracting the divergence at infinity, we may perform a partial fraction decomposition to get
\begin{equation*}
%\label{eq:cotransition_rational_function}
	H_{\lambda}(z) = z - \sum_{i=1}^{d-1}\frac{\wtilde{q}_{i}}{z-\y_{i}},\quad \wtilde{q}_{i} = -\frac{\prod_{j=1}^{d}(\y_{i}-\x_{j})}{\prod_{j=1;j\neq i}^{d-1}(\y_{i}-\y_{j})},\quad i=1,\dots, d-1.
\end{equation*}
%\cyan{[Check if this formula should be labelled.]}
There is no constant term appearing in the right-hand side because of (\ref{eq:maxmin_centered}),
and (\ref{eq:interlacing}) ensures that $\wtilde{q}_{i}>0$, $i=1,\dots, d-1$.
It can be shown that~\cite{kerov2000anisotropic}
\begin{equation*}
	\sum_{i=1}^{d-1}\wtilde{q}_{i} = |\lambda|,
\end{equation*}
thus, the sum $\sum_{i=1}^{d-1}\frac{1}{|\lambda|}\wtilde{q}_{i}\delta_{\y_{i}}$ is a probability measure on $\bR$
called the cotransition measure of $\lambda$.

%Notice that $y_{i}$'s are the contents of the boxes that are removable from $\lambda$.
Suppose that $\mu$ is another Young diagram such that $\lambda\searrow\mu$ with $c(\lambda/\mu)=\y_{i}$ with some $i$ and set
\begin{equation}
\label{eq:cotransition_normalized}
	p^{\downarrow}(\lambda,\mu):= \frac{1}{|\lambda|}\wtilde{q}_{i} = -\frac{1}{|\lambda|}\frac{\prod_{j=1}^{d}(\y_{i}-\x_{j})}{\prod_{j=1;j\neq i}^{d-1}(\y_{i}-\y_{j})}.
\end{equation}
%\cyan{[Labelling this equation might suffice in the future.]}
The probability measure on $\bY_{|\lambda|-1}$ defined by $p^{\downarrow}(\lambda,-)$ is also called the cotransition measure of $\lambda$.
Similarly to (\ref{eq:Cauchy_contents}), let us keep the following formula
\begin{equation}
\label{eq:cotrans_Cauchy_contents}
	H_{\lambda}(z) = z - |\lambda| \sum_{\mu;\lambda\searrow\mu} \frac{p^{\downarrow}(\lambda,\mu)}{z-c(\lambda/\mu)} .
\end{equation}

\subsection{Harmonic function}
We will need a harmonic function on the Young graph to define a planar algebra.
As the Young graph is considered an oriented graph, the following definition of a harmonic function is natural:
a harmonic function on $\bY$ is a non-negative function $f$ on $\bY$ such that
\begin{equation*}
    f(\mu) = \sum_{\lambda;\, \mu\nearrow\lambda}f(\lambda),\quad \mu\in\bY.
\end{equation*}
Classification of harmonic functions on $\bY$ has been completed in the context of the asymptotic representation theory of symmetric groups (see e.g.~\cite{kerov2004harmonic}).
Among others, the following example associated with the Plancherel measure will be important.

\begin{exam}
\label{exam:PlancherelHarmonic}
Let us define $f_{\Pl}\colon \bY \to \bR$ by
\begin{equation}
\label{eq:plancherel_harmonic}
	f_{\Pl}(\lambda) = \prod_{(i,j)\in\lambda}(\lambda_{i}+ \lambda'_{j}-i-j+1)^{-1},\quad \lambda\in \bY.
\end{equation}
Then, $f_{\Pl}$ is a harmonic function on $\bY$.
This harmonic function exhibits the following properties~\cite{kerov2000anisotropic}: for $\mu\nearrow \lambda$,
\begin{align}
\label{eq:ratio_f_transition}
	\frac{f_{\Pl}(\lambda)}{f_{\Pl}(\mu)}&= p^{\uparrow}(\mu,\lambda), \\
\label{eq:ratio_f_cotransition}
	\frac{f_{\Pl}(\mu)}{f_{\Pl}(\lambda)}&=|\lambda| p^{\downarrow}(\lambda,\mu).
\end{align}
\end{exam}

\begin{rem}
Perhaps, we need to explain how the harmonic function $f_{\Pl}$ is related to the Plancherel measure.
Most often, the Plancherel measure $\bM_{\Pl}$ is given by the formula
\begin{equation*}
	\bM_{\Pl}(\lambda) = \frac{(\dim \lambda)^{2}}{|\lambda|!},\quad \lambda\in \bY.
\end{equation*}
Here, $\dim$ is the dimension function given by the hook product formula
\begin{equation*}
%\label{eq:dimension_hook_length}
	\dim \lambda = |\lambda|!\cdot \prod_{(i,j)\in\lambda}(\lambda_{i}+ \lambda'_{j}-i-j+1)^{-1},\quad \lambda\in \bY.
\end{equation*}
%\cyan{[check if this is referred to.]}
The dimension of $\lambda$, $\dim \lambda$, is otherwise defined as the number of paths in $\bY$ from $\emptyset$ to $\lambda$, or equivalently, the number of standard tableaux of shape $\lambda$.
It is called dimension because $\dim \lambda$ is the dimension of the irreducible representation of the symmetric group that is labelled by $\lambda$.
It follows, then, from either the Burnside theorem or a direct calculation that the restriction of $\bM_{\Pl}$ to each $\bY_{n}$ gives a probability measure.

The cotransition probability $p^{\downarrow}(\lambda,\mu)$ is more usually defined by
\begin{equation}
\label{eq:cotrans_dim_ratio}
	p^{\downarrow}(\lambda,\mu) = \frac{\dim \mu}{\dim \lambda},\quad \lambda\searrow\mu.
\end{equation}
Then, it follows from the definition of the dimension (as the number of paths) that $p^{\downarrow}(\lambda,-)$ is indeed a probability measure under a fixed $\lambda$.

The Plancherel measure is said to form a coherent system in the sense that
\begin{equation*}
	\bM_{\Pl}(\mu) = \sum_{\lambda;\mu\nearrow\lambda}\bM_{\Pl}(\lambda)p^{\downarrow}(\lambda,\mu),\quad \mu\in \bY,
\end{equation*}
but as soon as we get a coherent system, we may divide it by the dimension function to get a harmonic function.
In fact, the harmonic function $f_{\Pl}$ is obtained by
\begin{equation*}
	f_{\Pl}(\lambda) = \frac{\bM_{\Pl}(\lambda)}{\dim \lambda},\quad \lambda\in \bY.
\end{equation*}

Now, we can see that (\ref{eq:cotrans_dim_ratio}) together with (\ref{eq:plancherel_harmonic}) recovers (\ref{eq:ratio_f_cotransition}), and (\ref{eq:ratio_f_transition}) gives
\begin{equation*}
	p^{\uparrow}(\mu,\lambda) = \frac{\dim \lambda}{|\lambda|\dim \mu},\quad \mu\nearrow \lambda
\end{equation*}
which is another common definition of the transition probabilities.
\end{rem}

\subsection{Moments and Boolean cumulants}
\label{sect:def_moments_cumulants}
Let $\lambda$ be a Young diagram.
The Cauchy transform $G_{\lambda}$ of $\sfm_{\lambda}$ is a generating function of the moments of $\sfm_{\lambda}$ as
\begin{equation}
\label{eq:def_moments}
	G_{\lambda}(z) = z^{-1} + \sum_{n=1}^{\infty} M_{n}(\lambda)z^{-n-1},\quad M_{n}(\lambda) = \int_{\bR}x^{n}\sfm_{\lambda}(dx),\, n\in \bZ_{>0}.
\end{equation}
Here, the right-hand side of the first equation is the expansion around infinity.

We have gotten numbers $M_{n}(\lambda)$, $n\in\bZ_{>0}$ for a given Young diagram $\lambda$.
Now, we switch the roles of $n$ and $\lambda$ to get a family of functions on $\bY$ labelled by $n\in \bZ_{>0}$.
For each $n\in \bZ_{>0}$, we will think of $M_{n}$ as a function on $\bY$ given by
\begin{equation*}
	M_{n}\colon \bY \to \bR;\quad \lambda\mapsto M_{n}(\lambda).
\end{equation*}

The Boolean cumulants $B_{n}(\lambda)$, $n\in\bZ_{>0}$ of $\lambda$ are generated by $H_{\lambda}$~\cite{SpeicherWoroudi1993}:
\begin{equation}
\label{eq:def_Boolean_cumulants}
	H_{\lambda}(z) = z - \sum_{n=1}^{\infty} B_{n}(\lambda)z^{-n+1},
\end{equation}
where the right-hand side is the Laurent expansion around infinity.
Again, we think of each $B_{n}$, $n\in\bZ_{>0}$ as a function on $\bY$ by
\begin{equation*}
	B_{n}\colon \bY \to \bR;\quad \lambda\mapsto B_{n}(\lambda).
\end{equation*}

\subsection{Normalized characters}
\label{subsect:normalizedCharacter}
For $n\in\bN$, let $\cP_{n}$ be the set of integer partitions of $n$.
Each element of $\cP_{n}$ can be presented as a non-increasing sequence of positive integers
$\pi = (\pi_{1},\dots, \pi_{l})$ sucht that $\sum_{i=1}^{l}\pi_{i}=n$.
We call the numbers $n$ and $l$ above the weight and length of the partition $\pi$ and write $|\pi|$ and $\ell (\pi)$, respectively.

The set $\cP_{n}$ can be naturally identified with the set $\bY_{n}$ of Young diagrams of the same weight,
but we prefer to distinguish them because of their different roles: $\cP_{n}$ labels the conjugacy classes in the symmetric group $S_{n}$, whereas $\bY_{n}$ labels the irreducible representations of the same group.
These two sets certainly have the same number of elements, but we consider the natural bijection $\cP_{n}\simeq \bY_{n}$ a mere coincidence.

%Characters can be defined in terms of representation theory of symmetric groups,
%but alternatively, we introduce them from symmetric functions.
%
Let $\lambda\in \bY_{n}$  with $n\in \bN$ be a Young diagram and $V^{\lambda}$ be the corresponding finite dimensional irreducible representation of $S_{n}$.
In particular, we have $\dim (V^{\lambda}) = \dim \lambda$.
The character of the representation is the function on $S_{n}$ defined by
\begin{align*}
	\chi^{\lambda} (\sigma) = \Tr_{V^{\lambda}}(\sigma),\quad \sigma\in S_{n}.
\end{align*}
The character only depends on conjugacy classes.
For $\pi\in\cP_{n}$, we set $\chi^{\lambda}_{\pi}=\chi^{\lambda}(\sigma)$,
where $\sigma$ lies in the conjugacy class labelled by $\pi$.

So far, partitions are variables of a character and Young diagrams are labels for characters.
We switch their roles now.
Let us introduce the set of all partitions $\cP = \bigcup_{n\geq 0}\cP_{n}$ with convention $\cP_{0}=\{\emptyset\}$.
What we call the normalized characters are functions on $\bY$ labelled by $\pi\in \cP$ defined as follows.
For $\pi\in \cP_{k}$ with $k\in\bN$, the normalized character $\Sigma_{\pi}$ is given by~\cite{Biane2003}
\begin{equation}
\label{eq:defNormalizedCharacter}
	\Sigma_{\pi}(\lambda) :=
	\begin{cases}
	(n\downharpoonright k)\frac{\chi^{\lambda}_{\pi\cup (1^{n-k})}}{\dim \lambda},& |\lambda|=n\geq k, \\
	0, & |\lambda|<k.
	\end{cases}
\end{equation}
Here, $(n\downharpoonright k) = n(n-1)\cdots (n-k+1)$ is the falling factorial
and $\pi\cup (1^{n-k})$ is the partition of $n$ obtained by adding $n-k$ copies of $1$ to $\pi$.

%The generalized Frobenius formula due to Rattan--{\'S}niady provides integral formulas for the normalized characters.
%\begin{thm}[Rattan--{\'S}niady]
%The normalized character $\Sigma_{\pi}$ for $\pi = (\pi_{1},\dots,\pi_{l})$ is computed by the following formula:
%\end{thm}

\section{Planar algebras}
\label{Sect:planarAlg}
In this section, we first introduce a general notion of planar algebra and next consider examples associated with the Young graph.
We follow the style of planar algebra due to Jones, but one feature of ours is that strings are oriented.
Thus, we might need to call our planar algebra an {\it oriented} planar algebra, but we just omit the word ``oriented'' as we do not need the unoriented version in this paper.
The idea of introducing orientation to planar algebras is by no means new.
It has been already mentioned in the original~\cite{jones2021planar},
and extensively used e.g. in \cite{bar2005khovanov}.
Nevertheless, we present a complete definition of it as details often differ from a literature to another.

\subsection{General setup}
\subsubsection{Planar tangles}
As we have sketched in Introduction, a planar algebra consists of vector spaces along with multi-linear maps labelled by planar tangles.
To make these notions precise, let us start by looking at planar tangles.

In the following, by a disk, we mean the image of the unit disk
\begin{equation*}
	\{(x,y)\in\bR^{2}|x^{2}+y^{2}<1\}
\end{equation*}
under a smooth embedding in $\bR^{2}$.
Although we occasionally draw a disk in an equiradial way, it does not mean that the disk must be isometrically embedded.

\begin{defn}
A planar tangle $T$ consists of the following data:
\begin{itemize}
\item 	A disk $D_{0}$ along with disjoint disks inside $D_{0}$: 
		\begin{equation*}
			D_{1},\dots, D_{N}\subset D_{0}.
		\end{equation*}
		We call $D_{0}$ the external disk, $D_{i}$, $i=1,\dots, N$ the internal disks, and write
		\begin{equation*}
			\cD_{T}=\{D_{i}: i=1,\dots, N\}
		\end{equation*}
		for the collection of the internal disks.
\item 	Each disk (external and internal) has marked boundary points with signs;
		each marked boundary point is denoted by $(p,\epsilon)$ where $p$ is a boundary point and $\epsilon \in \{+,-\}$.
		For each disk, the signs on its boundary need to be summed to be $0$,
		which, in particular, requires the number of boundary points on each disk to be even.
		We call a marked boundary point of the external disk an external boundary point and one of an internal disk an internal boundary point.
\item 	Non-intersectiong oriented strings in
		\begin{equation*}
		\Omega = D_{0}\bsl \bigcup_{i=1}^{N}D_{i}
		\end{equation*}
		obeying the following rules:
		\begin{itemize}
		\item Each string either forms a closed loop or connects two marked boundary points.
		\item A string connecting marked boundary points 
			\begin{itemize}
			\item starts either at an external marked boundary point with $\epsilon = -$ or an internal boundary point with $\epsilon = +$ and
			\item ends up either at an external boundary point with $\epsilon = +$ or an internal boundary point with $\epsilon = -$.
			\end{itemize}
		\item Each marked boundary point is connected to exactly one string.
		\end{itemize}
\item 	Each disk has a distinguished interval indicated by $*$ on its boundary separated by two of the marked boundary points.
\end{itemize}
We write $\bfT$ for the set of planar tangles.
\end{defn}

\begin{figure}
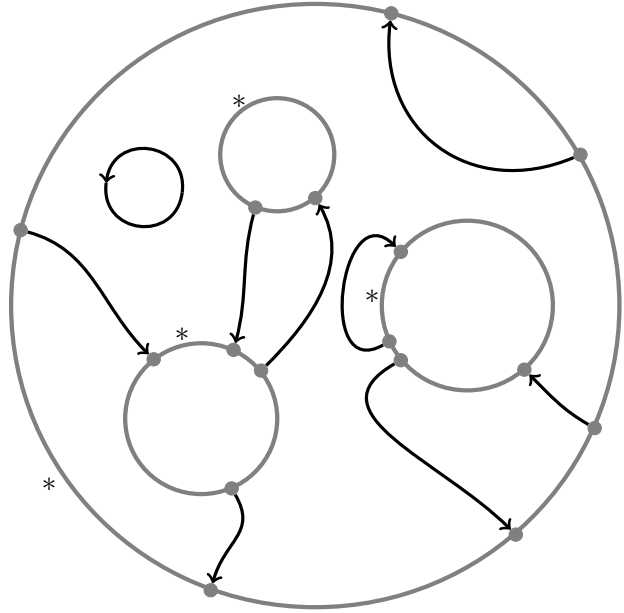

\ctikzfig{arcDiagram}
\caption{Planar tangle.}
\label{fig:example_oriented_planar_tangle}
\end{figure}

Figure~\ref{fig:example_oriented_planar_tangle} shows an example of a planar tangle.

\begin{rem}
It is common that a planar tangle is also equipped with shading of the regions separated by the strings
although, in this paper, we do no need shading.
%With shading, each disk must have even number of marked boundary points.
%In this paper, we do not need shading.
%However, for our relevant examples associated with the Young graph, the number of marked boundary points will be even.
\end{rem}

Given a planar tangle $T$ and a diffeomorphism $\Phi$, there is a natural way to send the whole structure of $T$ by $\Phi$.
In fact, we can send the disks, marked boundary points and strings by $\Phi$ and the additional data, namely, the signs of the marked boundary points and the distinguished boundary segments are also transferred consistently.
The outcome is a new planar tangle, and we write it as $\Phi (T)$.

\begin{figure}
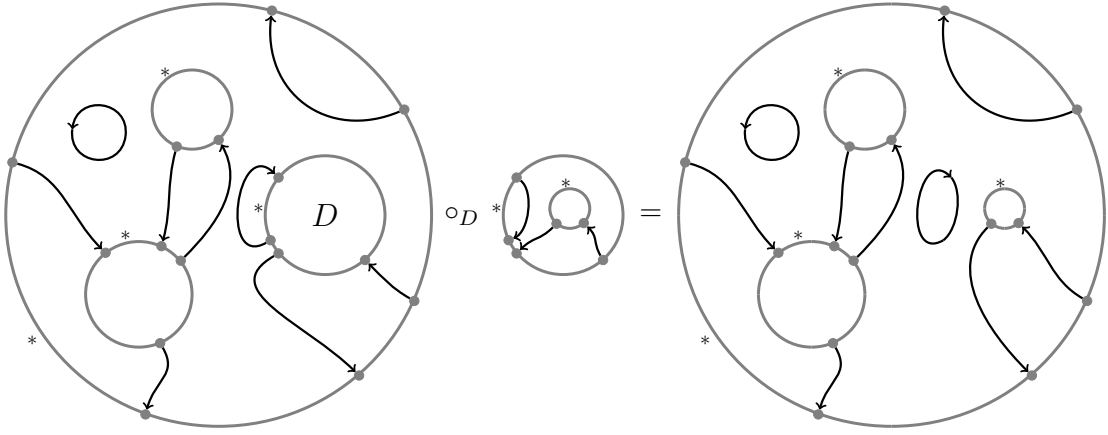

\begin{equation*}
\scalebox{.7}{\input{arcDiagramWithD.tikz}}\, \circ_{D}\, \scalebox{.7}{\input{arcDiagramSmall.tikz}} =
\scalebox{.7}{\input{arcDiagramComposed.tikz}}
\end{equation*}
\caption{Composition of planar tangles.}
\label{fig:composition_tangles}
\end{figure}

Let $T$ and $S$ be planar tangles and $D\in \cD_{T}$.
We say that $T$ is composable to $S$ over $D$ if the external disk of $S$ coincides with $D$ including the marked boundary points with signs and the distinguished boundary interval on it.
In this case, we can get a new planar tangle by filling $D$ in $T$ with $S$,
continuing the strings attached to the boundary of $D$, and forgetting the boundary of $D$.
The outcome will be denoted by $T\circ_{D} S$ and called the composition of $T$ and $S$ over $D$ (Figure~\ref{fig:composition_tangles}).

\subsubsection{Axioms for a planar algebra}
We will define a planar algebra as a family of vector spaces equipped with multi-linear maps labelled by planar tangles.
For this idea to make sense, the family of vector spaces needs to be labelled by all possible {\it boundary profiles} of planar tangles.
Since planar tangles have marked boundary points with signs, it is natural to suppose that the following set parametrizes the boundary profiles:
\begin{equation*}
	\bfE = \{\ul{\epsilon}=(\epsilon_{1},\dots, \epsilon_{n})| \epsilon_{i}\in \{+,-\},\, i=1,\dots, n,\, \textstyle{\sum_{i=1}^{n}}\epsilon_{i}=0\}.
\end{equation*}
Due to the last condition, the length $n$ must be even.
We agree that $\bfE$ contains a unique sequence of length $0$, which is denoted by $\ul{\emptyset}$.

Now, we can give a precise definition of our planar algebra.
A planar algebra is a family of vector spaces $(P_{\ul{\epsilon}})_{\ul{\epsilon}\in \bfE}$
equipped with a family of multi-linear maps in the following way.
Let $T$ be a planar tangle with $N$ number of internal disks.
Then, to each internal disk $D_{i}\in \cD_{T}$, $i=1,\dots, N$, we may assign a sequence of signs $\ul{\epsilon}^{(i)}\in\bfE$ by reading the signs of the marked boundary points on $D_{i}$ counterclockwise starting from the distinguished boundary interval.
We also get a sequence of signs $\ul{\epsilon}^{(0)}$ on the external disk in the same way.
Then, the multi-linear map associated with $T$ has the form
\begin{equation*}
    Z_{T}\colon P_{\ul{\epsilon}^{(1)}}\times \cdots \times P_{\ul{\epsilon}^{(N)}} \to P_{\ul{\epsilon}^{(0)}}.
\end{equation*}
In the case that $T$ has no internal disk, we understand the left-hand side as the scalar.

\begin{defn}
Let $P=(P_{\ul{\epsilon}})_{\ul{\epsilon}\in\bfE}$ be a family of vector spaces and $Z=(Z_{T})_{T\in\bfT}$ be a family of multi-linear maps of the above type.
The pair $(P,Z)$ forms a planar algebra if the following properties are satisfied:
\begin{enumerate}
    \item Isotopy invariance: under any orientation preserving diffeomorphism $\Phi$ of $\bR^{2}$, we have 
    \begin{equation*}
        Z_{T}=Z_{\Phi (T)}
    \end{equation*}
    for all planar tangles $T\in\bfT$.
    \item Consistency under composition: if $T$ is composable to $S$ over $D\in\cD_{T}$,
    \begin{equation*}
        Z_{T\circ_{D} S}= Z_{T}\circ \iota_{D}(Z_{S}).
    \end{equation*}
    Here, $\iota_{D}(Z_{S})$ makes $Z_{S}$ act on the vector spaces corresponding to the internal disks of $S$
    and returns the outcome to the vector space corresponding to $D$.
\end{enumerate}
\end{defn}

We can always isotope a planar tangle so that each disk has
\begin{itemize}
\item 	the shape of a rectangle with smooth corners,
\item 	marked boundary points on its {\it top edge},
\item 	and its distinguished boundary segments containing the left, bottom and right edges.
\end{itemize}

\begin{defn}
We say that a planar tangle satisfying the above three conditions is of a standard shape.
\end{defn}

Note that, for a planar tangle of a standard shape,
no information is lost if we replace the smooth corners of rectangles with sharp ones
while it is not realized by a diffeomorphism of $\bR^{2}$.
The right hand side of the following shows a typical oriented planar tangle under the above identifications:
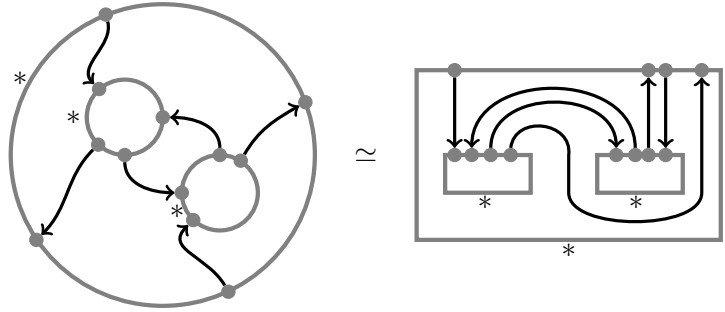
\begin{equation*}
\input{arcDDemoDisk.tikz} \quad \simeq \quad  \input{arcDDemoRect.tikz}\quad .
\end{equation*}
We will present planar tangles of standard shapes in this way mainly due to the limitation of the author's drawing skills.

\subsection{Example associated with the Young graph}
We see an example of constructing a planar algebra that is associated with the Young graph.
The construction here is a natural extension of the one of an (unoriented) planar algebra for a bipartite graph~\cite{jones2000planar}
to the case of an oriented graph.
We fix a function $f\colon \bY \to \bR_{>0}$ that will be used to define the multi-linear maps for the planar algebra.

\begin{rem}
We assumed that the function $f$ is non-vanishing for simplicity.
This condition can be relaxed so that the values of $f$ are non-negative
by using the support of $f$ as our oriented graph instead of the whole $\bY$.
\end{rem}

\subsubsection{Vector spaces}
The first ingredient for a planar algebra is a family of vector spaces $P=(P_{\ul{\epsilon}})_{\ul{\epsilon}\in\bfE}$.
For a sequence of signs $\ul{\epsilon}\in\bfE$,
we set $P_{\ul{\epsilon}}$ to be the vector space of functions of loops on $\bY$ of signature $\ul{\epsilon}$ (see Section~\ref{subsect:YoungGraph}).
Note that, in particular, $P_{\ul{\emptyset}}$ consists of functions on $\bY$, thus $P_{\ul{\emptyset}}\simeq \bC [\bY]$.

%anchored loops of length $k$ in $\bY$ of the form
%\begin{equation*}
%    \lambda^{(0)}\xrightarrow{\epsilon_{1}} \lambda^{(1)} \xrightarrow{\epsilon_{2}} \lambda^{(2)} \xrightarrow{\epsilon_{3}}\cdots \xrightarrow{\epsilon_{k}}\lambda^{(k)}= \lambda^{(0)},
%\end{equation*}
%where we used the convention of arrows as
%\begin{equation*}
%    \xrightarrow{\epsilon} \, =
%    \begin{cases}
%        \searrow\, , & \epsilon = +, \\
%        \nearrow\, , & \epsilon = -.
%    \end{cases}
%\end{equation*}
%Since there is no loop of odd length, $P_{\ul{\epsilon}}=0$ if $k$ is odd.

\subsubsection{Spin states and weights}
The other member of a planar algebra is a family of multi-linear maps $Z = (Z_{T})_{T\in \bfT}$.
We construct them in terms of state-sum as we will see below.

Let $T$ be a planar tangle.
A spin state on $T$ is a filling of the regions of $T$ separated by the strings with Young diagrams obeying the following rule.
When standing on each string along the orientation, let $\lambda$ and $\mu$ be the Young diagrams on the left-hand side and on the right-hand side, respectively:
\[
\begin{tikzpicture}
	\begin{pgfonlayer}{nodelayer}
		\node [style=none] (0) at (0, -1) {};
		\node [style=none] (2) at (0, 1) {};
		\node [style=none] (3) at (-1, 0) {$\lambda$};
		\node [style=none] (4) at (1, 0) {$\mu$};
	\end{pgfonlayer}
	\begin{pgfonlayer}{edgelayer}
		\draw [style=verythickArrow] (0.center) to (2.center);
	\end{pgfonlayer}
\end{tikzpicture}
\quad .
\]
Then, we must have $\mu\nearrow\lambda$.

Given a spin state on a planar tangle, the boundary profile of the spin state is the sequence of Young diagrams that fill the regions connected to the external boundary read counterclockwise starting from the distinguished interval.
Of course, the boundary profile forms a loop
\begin{equation*}
    \ul{\lambda}=(\lambda^{(0)}\xrightarrow{\epsilon_{1}} \lambda^{(1)} \xrightarrow{\epsilon_{2}} \lambda^{(2)} \xrightarrow{\epsilon_{3}}\cdots \xrightarrow{\epsilon_{k}}\lambda^{(n)}= \lambda^{(0)})
\end{equation*}
in $\bY$, where $\ul{\epsilon}\in \bfE$ is the sequence of signs on the external disk.
The collection of spin states on $T$ with a boundary profile $\ul{\lambda}$ is denoted by $\Conf_{\ul{\lambda}}(T)$ (standing for {\it configuration}).
Note that $\Conf_{\ul{\lambda}}(T)$ is a finite set.

For example, the filling
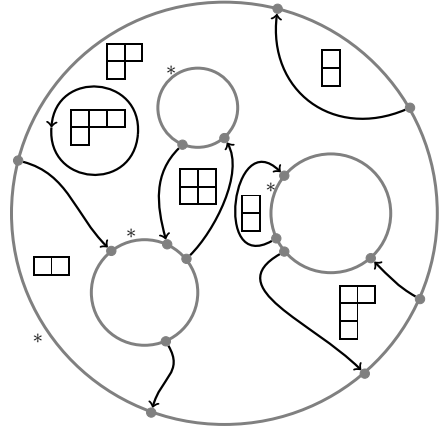
\[
\scalebox{.7}{\input{spinStateExample.tikz}}
\]
is a spin state with the boundary profile
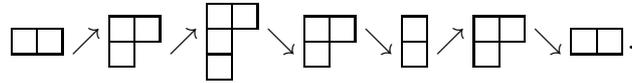
\begin{equation*}
	\hackcenter{
	\begin{ytableau}
	\, & \,
	\end{ytableau}
	} \nearrow 
	\hackcenter{
	\begin{ytableau}
	\, & \, \\
	\,
	\end{ytableau}
	} \nearrow
	\hackcenter{
	\begin{ytableau}
	\,  &\, \\
	\, \\
	\,
	\end{ytableau}
	} \searrow 
	\hackcenter{
	\begin{ytableau}
	\, & \, \\
	\,
	\end{ytableau}
	} \searrow
	\hackcenter{
	\begin{ytableau}
	\,  \\
	\,
	\end{ytableau}
	} \nearrow
	\hackcenter{
	\begin{ytableau}
	\, & \, \\
	\,
	\end{ytableau}
	} \searrow
	\hackcenter{
	\begin{ytableau}
	\, & \,
	\end{ytableau}
	}\, .
\end{equation*}

We next define the weight of a spin state $\sigma$ on a planar tangle $T$.
For that, we first isotope $T$ to a planar tangle of a standard shape in which
the critical points of the strings (in terms of the height against the vertical direction) are either maximal or minimal.
Then, the spin state $\sigma$ on $T$ naturally induces a spin state on the resulting oriented planar tangle after isotopy.
For each critical point of a string, we find the Young diagram $\lambda_{\mathrm{convex}}$ on the convex side and $\lambda_{\mathrm{concave}}$ on the concave side that are assigned by the spin state $\sigma$:
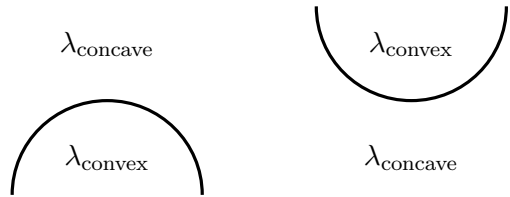
\begin{align*}
	\input{convexConcave.tikz}
\end{align*}
In this picture, we dropped the orientation of the string since the weights are defined independently of the orientations.
We assign to this critical point the local weight
\begin{equation*}
	\sqrt{\frac{f(\lambda_{\mathrm{convex}})}{f(\lambda_{\mathrm{concave}})}}.
\end{equation*}
Here, $f$ is the positive real valued function on $\bY$ that we have fixed in the beginning.
Then, the weight of the spin state $\sigma$ is the product of those local weights over all critical points:
\begin{equation*}
	\wt_{f} (\sigma) = \prod_{\text{critical points}}\sqrt{\frac{f(\lambda_{\mathrm{convex}})}{f(\lambda_{\mathrm{concave}})}}.
\end{equation*}

\subsubsection{Multilinear maps}
Now we are ready to define multilinear maps labelled by planar tangles.
Let $T$ be a planar tangle with $N$ internal disks: $\cD_{T}=\{D_{1},\dots,D_{N}\}$.
For each $i=1,\dots, N$, the sequence of signs on the boundary of $D_{i}$ is denoted by $\ul{\epsilon}^{(i)}\in \bfE$.
Also, we set $\ul{\epsilon}^{(0)}$ for the sequence of signs on the external disk $D_{0}$.

For $i=1,\dots, N$, let $h_{i}\in P_{\ul{\epsilon}^{(i)}}$ be a function attached to the internal disk $D_{i}$.
Then, we can {\it evaluate} $h_{i}$ on a spin state $\sigma$ on $T$ as follows.
When we read the Young diagrams assigned by $\sigma$ to the regions connected to the boundary of $D_{i}$ counterclockwise,
we get a loop in $\bY$ of signature $\ul{\epsilon}^{(i)}$.
Then, $h_{i}$ can be evaluated on the loop and the result can be written as $h_{i}(\sigma)$.

We define a multilinear map
\begin{equation*}
	Z^{f}_{T}\colon P_{\ul{\epsilon}^{(1)}}\times \cdots \times P_{\ul{\epsilon}^{(N)}}\to P_{\ul{\epsilon}^{(0)}}
\end{equation*}
by the formula
\begin{equation}
\label{eq:defMultilinearMap}
	Z^{f}_{T}(h_{1},\dots, h_{N})(\ul{\lambda})=\sum_{\sigma\in \Conf_{\ul{\lambda}}(T)}\Bigl(\prod_{i=1}^{N}h_{i}(\sigma)\Bigr) \wt_{f} (\sigma),
\end{equation}
where $h_{i}\in P_{\ul{\epsilon}^{(i)}}$, $i=1,\dots, N$ and $\ul{\lambda}$ is a loop in $\bY$ of signature $\ul{\epsilon}^{(0)}$.
In the case that $T$ has no internal disk, (\ref{eq:defMultilinearMap}) requires an interpretation.
In such a case, the product over empty set is understood as $\bC$ and $Z^{f}_{T}$ is determined by the value at $1$.
The right-hand of (\ref{eq:defMultilinearMap}) still makes sense as the definition of $Z^{f}_{T}(1)(\ul{\lambda})$ by understanding the product of $h_{i}$ as unity.

\begin{prop}
The pair of a family of vector spaces $P=(P_{\ul{\epsilon}})_{\ul{\epsilon}\in \bfE}$
and one of multi-linear maps $Z^{f}=(Z^{f}_{T})_{T\in\bfT}$ forms a planar algebra.
\end{prop}
\begin{proof}
We need to prove the isotopy invariance and the consistency under composition.

Let us start by the isotopy invariance.
It suffices to show that the definition of the weight of a spin state does not depend on a way to isotope a planar tangle to one of a standard shape.
For two planar tangles of standard shapes, if they are isotopic, the isotopy between them is a composition of bending a string and rotating a box by 360$^{\circ}$.
It is rather clear that the weight of a spin state is invariant under bending a string;
the contributions from the two critical points cancel with each other:
\begin{align*}
	\input{bendingString.tikz}\quad .
\end{align*}
The invariance under rotation also follows from cancellation;
suppose that, given a spin state, the loop in $\bY$ along the boundary of the rotated box is
\begin{equation*}
	\lambda^{(0)} \to \lambda^{(1)}\to \lambda^{(2)} \to \cdots \to \lambda^{(n)}=\lambda^{(0)},
\end{equation*}
where we omit the signs on the arrows as they do not matter here.
Then, after a clockwise rotation, the strings connected to the box contribute to the weight
\begin{equation*}
	\frac{f(\lambda^{(0)})}{f(\lambda^{(1)})}\frac{f(\lambda^{(1)})}{f(\lambda^{(2)})}\cdots \frac{f(\lambda^{(n-1)})}{f(\lambda^{(n)})} = \frac{f(\lambda^{(0)})}{f(\lambda^{(n)})} = 1.
\end{equation*}
Thus, the weight is the same as in the case without rotation:
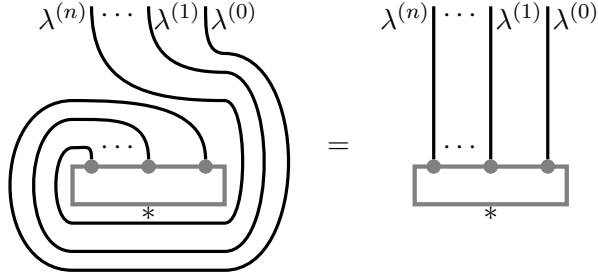
\begin{align*}
	\input{rotationBox.tikz}\quad .
\end{align*}
The invariance under a counterclockwise rotation is shown by the same way.

Let us move on to the consistency under composition.
Let $T$ be a planar tangle with $N$ internal disks: $\cD_{T} = \{D_{1},\dots, D_{N}\}$
and $S$ be one with $M$ internal disks $\cD_{S} = \{D_{N+1},\dots, D_{N+M}\}$.
We assume that $T$ is composable to $S$ over $D_{N}$.
Here, choosing the internal disk $D_{N}$ does not lose generality.
For each $i=1,\dots, N+M$, we write $\ul{\epsilon}^{(i)}\in \bfE$ for the sequence of signs on the boundary of $D_{i}$,
and $\ul{\epsilon}^{(0)}$ for that on the external disk of $T$.
The goal is to show the identity
\begin{align*}
	&\, Z^{f}_{T\circ_{D_{N}}S}(h_{1},\dots,h_{N-1},h_{N+1},\dots,h_{N+M}) \\
	=&\, Z^{f}_{T}(h_{1},\dots, h_{N-1},Z_{S}(h_{N+1},\dots, h_{N+M}))
\end{align*}
in $P_{\ul{\epsilon}^{(0)}}$ for all $h_{i}\in P_{\ul{\epsilon}^{(i)}}$, $i=1,\dots,N-1,N+1,\dots,N+M$.

Now, we take a loop $\ul{\lambda}$ in $\bY$ of signature $\ul{\epsilon}^{(0)}$.
Given a spin state $\sigma\in \Conf_{\ul{\lambda}}(T\circ_{D_{N}}S)$, we can get spin states on $T$ and $S$ by restriction. Let us write those restrictions as $\sigma|_{T}$ and $\sigma|_{S}$, respectively.
For each loop $\ul{\mu}$ in $\bY$ of signature $\ul{\epsilon}^{(N)}$, we set
\begin{equation*}
	\Conf_{\ul{\lambda},(\ul{\mu},N)}(T) = \{\sigma \in \Conf_{\ul{\lambda}}(T)|\text{reading $\ul{\mu}$ on the boundary of $D_{N}$}\}.
\end{equation*}
Then, we have a bijection
\begin{equation*}
	\Conf_{\ul{\lambda}}(T\circ_{D_{N}}S) \simeq \bigsqcup_{\ul{\mu}}\Conf_{\ul{\lambda},(\ul{\mu},N)}(T)\times \Conf_{\ul{\mu}}(S);\quad \sigma \mapsto (\sigma|_{T},\sigma|_{S}),
\end{equation*}
where $\ul{\mu}$ runs over all loops in $\bY$ of signature $\ul{\epsilon}^{(N)}$.
Since, the weight of a spin state is given by a product of local contributions, we can see that
\begin{equation*}
	\wt_{f} (\sigma) = \wt_{f} (\sigma|_{T}) \wt_{f} (\sigma|_{S}),\quad \sigma\in \Conf_{\ul{\lambda}}(T\circ_{D_{N}}S).
\end{equation*}

The desired identity is now straightforward:
for a loop $\ul{\lambda}$ in $\bY$ of signature $\ul{\epsilon}^{(0)}$,
we have
\begin{align*}
	&\,Z^{f}_{T\circ_{D_{N}}S}(h_{1},\dots,h_{N-1},h_{N+1},\dots,h_{N+M})(\ul{\lambda}) \\
	=&\, \sum_{\sigma\in\Conf_{\ul{\lambda}(T\circ_{D_{N}}S)}} \Bigl(\prod_{i=1}^{N-1}h_{i}(\sigma)\Bigr)\Bigl(\prod_{i=N+1}^{N+M}h_{i}(\sigma)\Bigr)\wt_{f} (\sigma) \\
	=&\, \sum_{\ul{\mu}}\sum_{\sigma_{1}\in\Conf_{\ul{\lambda},(\ul{\mu},N)}(T)}\Bigl(\prod_{i=1}^{N-1}h_{i}(\sigma_{1})\Bigr)\wt_{f} (\sigma_{1}) \sum_{\sigma_{2}\in\Conf_{\ul{\mu}}(S)}\Bigl(\prod_{i=N+1}^{N+M}h_{i}(\sigma_{2})\Bigr)\wt_{f} (\sigma_{2}) \\
	=&\, \sum_{\ul{\mu}}\sum_{\sigma_{1}\in\Conf_{\ul{\lambda},(\ul{\mu},N)}(T)}\Bigl(\prod_{i=1}^{N-1}h_{i}(\sigma_{1})\Bigr)Z^{f}_{S}(h_{N+1},\dots,h_{N+M})(\ul{\mu})\wt_{f} (\sigma_{1}) \\
	=&\, \sum_{\sigma\in\Conf_{\ul{\lambda}}(T)}\Bigl(\prod_{i=1}^{N-1}h_{i}(\sigma)\Bigr)Z^{f}_{S}(h_{N+1},\dots,h_{N+M})(\sigma)\wt_{f} (\sigma) \\
	=&\, Z^{f}_{T}(h_{1},\dots,h_{N-1},Z^{f}_{S}(h_{N+1},\dots,h_{N+M}))(\ul{\lambda}).
\end{align*}
Thus, the proof is complete.
\end{proof}

\subsubsection{Case with no internal disk}
\label{Sect:tangleWOIntDisk}
As we have already mentioned, if a planar tangle $T$ has no internal disk, $Z^{f}_{T}$ is a linear map of the form
\begin{equation*}
	Z^{f}_{T}\colon \bC \to P_{\ul{\epsilon}^{(0)}}.
\end{equation*}
Thus, we may identify $Z^{f}_{T}$ with $Z^{f}_{T}(1)\in P_{\ul{\epsilon}^{(0)}}$ that is defined by the formula (\ref{eq:defMultilinearMap}) along with the interpretation following it.
This allows us to consider a planar tangle without an internal disk a function of loops in $\bY$.
Just as an example, the tangle
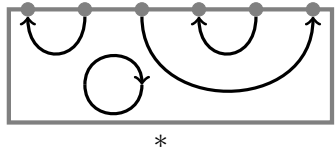
\begin{align*}
	\input{tangleFunction.tikz}
\end{align*}
is considered a function of loops of signature $(+,-,+,-,-,+)$ whose values are given by
\begin{align*}
	\input{tangleFunctionEval.tikz}\quad =&\, \sum_{\mu;\mu\nearrow \lambda^{(0)}} \input{tangleFunctionEvalSum.tikz} \\[1em]
	=&\, \sum_{\mu;\mu\nearrow \lambda^{(0)}}\sqrt{\frac{f(\lambda^{(2)})}{f(\lambda^{(0)})}}\sqrt{\frac{f(\lambda^{(5)})}{f(\lambda^{(0)})}}\frac{f(\mu)}{f(\lambda^{(0)})}\delta_{\lambda^{(1)},\lambda^{(3)}}\delta_{\lambda^{(4)},\lambda^{(0)}}.
\end{align*}

To obtain the planar algebra $(P,Z^{f})$, it sufficed to assume that the function $f$ is positive real valued.
Here, we add a further assumption that $f$ is harmonic.
With this assumption, we may recover the second relation in (\ref{eq:KHrelation3}) in the setting of a planar algebra.
This is analogous to the loop removing in the case of a bipartite graph~\cite{jones2000planar}.
\begin{prop}
\label{prop:PALeftCircle}
Let $f$ be a harmonic function on the Young graph and $(P,Z^{f})$ be the corresponding planar algebra.
As an element of $P_{\ul{\emptyset}}$,
\begin{align*}
	\input{leftCircle.tikz} = 1
\end{align*}
is a constant function on $\bY$.
\end{prop}
\begin{proof}
For an element of $P_{\ul{\emptyset}}$, we may omit the external disk.
With this convention, we have
\begin{align*}
	\input{leftCircleEval.tikz}\quad =&\, \sum_{\mu;\lambda\nearrow\mu}\input{leftCircleEvalSum.tikz}\quad 
	=\sum_{\mu;\lambda\nearrow\mu}\frac{f(\mu)}{f(\lambda)}=1.
\end{align*}
The equality follows from the harmonicity of $f$.
\end{proof}

\subsubsection{Comments on convention}
\label{sect:planar_alg_convention}
In the above arguments, we have implicitly applied the following convention:
when we write or draw Young diagrams in a planar tangle or part of it,
it refers to the evaluation of the planar tangle or part of it with those Young diagrams
according to (\ref{eq:defMultilinearMap}).
This convention will be used for the rest of the paper.

In the proof of Proposition~\ref{prop:PALeftCircle}, we omitted the external disk for an element of $P_{\ul{\emptyset}}$.
This is also legitimate since no information is lost by doing so.
We will adopt this as well in the sequel.

\section{Example: Plancherel case}
\label{Sect:plancherelCase}
In this section, we take the harmonic function $f_{\Pl}$ on $\bY$ pertaining to the Plancherel measure introduced in Example~\ref{exam:PlancherelHarmonic} and study the associated oriented planar algebra.
We first introduce an element in the planar algebra corresponding to a crossing of two strings,
and then present local relations satisfied by planar algebra elements that recover (\ref{eq:KHrelation1})--(\ref{eq:KHrelation3}).
We also give planar algebra elements that are identified with moments, Boolean cumulants, and normalized characters of Young diagrams.

\subsection{Crossing element}
We introduce a special element in $P_{(-,-,+,+)}$ that we will call the crossing element,
but first start by defining its building blocks
\begin{equation*}
	t_{\id}=\input{crossId.tikz}\quad ,\qquad t_{\ex}=\input{crossExchange.tikz}\quad .
\end{equation*}
They both are elements of $P_{(-,-,+,+)}$ and to be evaluated on loops of signature $(-,-,+,+)$:
\begin{equation*}
	\ul{\lambda} = (\lambda^{(0)} \searrow \lambda^{(1)} \searrow \lambda^{(2)} \nearrow \lambda^{(3)}\nearrow \lambda^{(4)}= \lambda^{(0)}).
\end{equation*}
For such loops, we define $t_{\id}$ by
\begin{align*}
	t_{\id} (\ul{\lambda})
		=&\, \input{crossIdEval.tikz} \\
	 :=&\, (c(\lambda^{(0)}/\lambda^{(1)})-c(\lambda^{(1)}/\lambda^{(2)}))^{-1}\sqrt{\frac{f_{\Pl}(\lambda^{(2)})}{f_{\Pl}(\lambda^{(0)})}} \delta_{\lambda^{(1)},\lambda^{(3)}}.
\end{align*}
Note that we have $c(\lambda^{(0)}/\lambda^{(1)})\neq c(\lambda^{(1)}/\lambda^{(2)})$, so the right-hand side is well-defined.
The subscript $\id$ stands for ``identity'' as the paths $\lambda^{(2)}\nearrow \lambda^{(1)} \nearrow \lambda^{(0)}$ and $\lambda^{(2)}\nearrow \lambda^{(3)} \nearrow \lambda^{(0)}$ are identical as long as $t_{\id}(\ul{\lambda})\neq 0$.

The other one $t_{\ex}$ is defined on loops $\ul{\lambda}$ of the same form by
\begin{align*}
	t_{\ex}(\ul{\lambda})
		=&\, \input{crossExchangeEval.tikz}\\
		:=&\,\sqrt{1-(c(\lambda^{(0)}/\lambda^{(1)})-c(\lambda^{(1)}/\lambda^{(2)}))^{-2}}\sqrt{\frac{f_{\Pl}(\lambda^{(2)})}{f_{\Pl}(\lambda^{(0)})}}(1-\delta_{\lambda^{(1)},\lambda^{(3)}}).
\end{align*}
The subscript $\ex$ stands for ``exchange'' this time because of the following reason:
for a loop $\ul{\lambda}$ such that $t_{\ex}(\ul{\lambda})\neq 0$,
when we compare the paths $\lambda^{(2)}\nearrow \lambda^{(1)} \nearrow \lambda^{(0)}$ and $\lambda^{(2)}\nearrow \lambda^{(3)} \nearrow \lambda^{(0)}$, they differ by the order of adding the two boxes.

\begin{rem}
It is a good mental picture that each string is {\it carrying} a box.
In fact, for any loop of signature $(-,-,+,+)$, we have
\begin{align*}
	\input{crossIdEval.tikz}\, \neq 0 &\Leftrightarrow \lambda^{(0)}/\lambda^{(1)}=\lambda^{(0)}/\lambda^{(3)} \text{ and } \lambda^{(1)}/\lambda^{(2)}=\lambda^{(3)}/\lambda^{(2)}, \\
	\input{crossExchangeEval.tikz}\, \neq 0 &\Leftrightarrow \lambda^{(0)}/\lambda^{(1)}=\lambda^{(3)}/\lambda^{(2)} \text{ and } \lambda^{(1)}/\lambda^{(2)}=\lambda^{(0)}/\lambda^{(3)}.
\end{align*}
\end{rem}

Finally, we define the crossing element $t\in P_{(-,-,+,+)}$ as the sum of $t_{\id}$ and $t_{\ex}$:
\begin{align*}
	t =\input{crossing.tikz}
	:=t_{\id}+t_{\ex} = \input{crossId.tikz} + \input{crossExchange.tikz}\, .
\end{align*}

We will use the graphical expression of the above elements in the example illustrated below.
Suppose that a tangle $T$ is
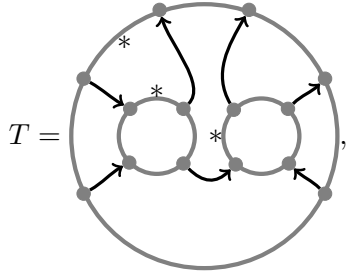
\begin{equation*}
	T = \input{tangleExamCrossDemo.tikz},
\end{equation*}
the multi-linear map $Z_{T}$ can be evaluated at $(t,t)\in (P_{(-,-,+,+)})^{2}$ and the result can be drawn as
\begin{equation*}
	Z_{T}(t,t) = \input{crossEmbed.tikz} \in P_{(-,-,-,+,+,+)},
\end{equation*}
or we can even omit the internal disks of $T$ and isotope the diagram into
\begin{equation*}
	Z_{T}(t,t) = \input{crossEmbedSimple.tikz}
\end{equation*}
as long as we remember the location of $*$ relative to each crossing.
In the sequel, we will omit the external disk of $t$ when we depict the result of evaluating it by a multi-linear map $Z_{T}$.
We will similarly draw $t_{\id}$ and $t_{\ex}$ embedded in a larger planar tangle as well,
but for them we cannot omit their external disks.

\begin{rem}
We shall make a couple of comments on the above definitions.
\begin{enumerate}
\item Apart from the part $\sqrt{\frac{f_{\Pl}(\lambda^{(2)})}{f_{\Pl}(\lambda^{(0)})}}$ involving the harmonic function,
the rest in the definition of $t$ reads the matrix elements of a simple transpose of a symmetric group acting on an irreducible representation (see Appendix~\ref{app:reptheorSG}).
We will not need this fact for the most part except in Section~\ref{sect:normalizedCharacters}.

\item The reason for having the part $\sqrt{\frac{f_{\Pl}(\lambda^{(2)})}{f_{\Pl}(\lambda^{(0)})}}$ is
merely because of our definition of the multi-linear maps $Z_{T}$.
In fact, $t$ will be evaluated in a rectangle with four boundary points on the top edge,
so will be {\it naturally isotoped} into
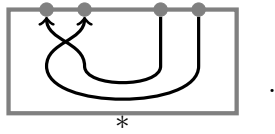
\begin{equation*}
	\input{crossingRect.tikz}\quad .
\end{equation*}
In this picture, we find two critical points on the strings that give rise to the harmonic functions.
\end{enumerate}
\end{rem}

We have a particular harmonic function $f_{\Pl}$ in the definition of $t$.
Just for defining $t$, we can replace it with another harmonic function,
but it will turn out that several interesting properties hold under this choice.

\subsection{Local relations}
We are now ready to see local relations exhibited by our planar algebra.
The relations we are about to show are those defining relations (\ref{eq:KHrelation1})--(\ref{eq:KHrelation3}) for the Khovanov Heisenberg category.
As we have mentioned in Introduction, in principle, the proof can be done by reformulating our planar algebra in the language of representations of symmetric groups.
We would rather not do so, and present our proof independent of the representation theory to feature the planar algebraic perspective.

\begin{thm}
\label{thm:local_relations}
Let $(P,Z^{f_{\Pl}})$ be the planar algebra associated with the harmonic function $f_{\Pl}$.
We have the following relations among elements in the planar algebra.
\begin{description}
\item[R1] %\label{item:leftTurn}
	In $P_{(-,+)}$, we have
	\begin{equation*}
		\input{leftTurn.tikz}\quad = 0.
	\end{equation*}
\item[R2-1] %\label{item:indInd}
	In $P_{(-,-,+,+)}$, we have
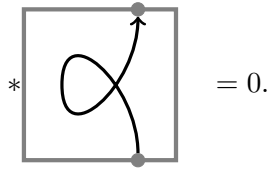
	\begin{equation*}
		\input{indIndCrossings.tikz}\quad = \quad \input{indIndStraight.tikz}\quad .
	\end{equation*}
\item[R2-2] %\label{item:indRes}
	In $P_{(-,+,-,+)}$, we have
	\begin{equation*}
		\input{indResCrosses.tikz} \quad = \quad \input{indResStraight.tikz}\quad .
	\end{equation*}
\item[R2-3] %\label{item:resInd}
	In $P_{(+,-,+,-)}$, we have
	\begin{equation*}
		\input{resIndCrosses.tikz} \quad =\quad \input{resIndStraight.tikz} \quad - \quad \input{resIndCups.tikz}\quad .
	\end{equation*}
\item[R3] %\label{item:YBE}
	In $P_{(-,-,-,+,+,+)}$, we have
	\begin{equation*}
		\input{yangBaxterLeft.tikz}\quad = \quad \input{yangBaxterRight.tikz}\quad .
	\end{equation*}
\end{description}
\end{thm}

Recall Section~\ref{Sect:tangleWOIntDisk} for the way to understand planar tangles with no internal disk
as elements in the planar algebra.

We will move on to a new subsection for the proof of Theorem~\ref{thm:local_relations}.
Notice that Theorem~\ref{thm:local_relations} combined with Proposition~\ref{prop:PALeftCircle}
recovers the defining relations for the Khovanov Heisenberg category (\ref{eq:KHrelation1})--(\ref{eq:KHrelation3}).

\begin{rem}
In the local relations in Theorem~\ref{thm:local_relations}, we have fixed distinguished intervals in a uniform way.
Note that we can change the location of the distinguished intervals simultaneous by applying the planar tangle
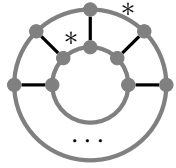
\begin{align*}
	\input{rotation.tikz}
\end{align*}
on both sides.
\end{rem}

\subsection{Proof of Theorem~\ref{thm:local_relations}}
\subsubsection{Preliminaries}
It will be convenient to prepare a few computational lemmas.

\begin{lem}
\label{lem:second_order_pole_transition}
Suppose that $\mu\nearrow\lambda$ in $\bY$.
We have
\begin{align*}
	\sum_{\nu;\lambda\nearrow \nu}\frac{p^{\uparrow}(\lambda,\nu)}{(c(\nu/\lambda)-c(\lambda/\mu))^{2}}
	= \frac{1}{|\lambda|p^{\downarrow}(\lambda,\mu)}.
\end{align*}
\end{lem}
\begin{proof}
We differentiate both sides of (\ref{eq:Cauchy_contents}) in the variable $z$ to get
\begin{equation*}
	G^{'}_{\lambda} (z) = -\sum_{\nu;\lambda \nearrow \nu}\frac{p^{\uparrow} (\lambda,\nu)}{(z-c(\nu/\lambda))^{2}}.
\end{equation*}
We also have an explicit formula for the Cauchy transform $G_{\lambda}$ from (\ref{eq:def_of_G_rational}), thus we obtain
\begin{align*}
	&\,\sum_{\nu;\lambda \nearrow \nu}\frac{p^{\uparrow} (\lambda,\nu)}{(z-c(\nu/\lambda))^{2}} \\
	=&\, -\sum_{i=1}^{d-1}\frac{\prod_{j=1;j\neq i}^{d-1}(z-\y_{j})}{\prod_{j=1}^{d}(z-\x_{j})} + (\text{terms vanishing at $\y_{i}$, $i=1,\dots, d-1$}).
\end{align*}
Here, $\x_{1},\dots, \x_{d}$ and $\y_{1},\dots, \y_{d-1}$ are the local minima and maxima of the profile of $\lambda$.
Now, we substitute $z = c(\lambda /\mu)=\y_{i}$ with some $i=1,\dots,d-1$ and get
\begin{align*}
	\sum_{\nu;\lambda\nearrow \nu}\frac{p^{\uparrow}(\lambda,\nu)}{(c(\nu/\lambda)-c(\lambda/\mu))^{2}}
	= -\frac{\prod_{j=1;j\neq i}^{d-1}(\y_{i}-\y_{j})}{\prod_{j=1}^{d}(\y_{i}-\x_{j})}.
\end{align*}
When we compare the right-hand side with (\ref{eq:cotransition_normalized}),
we find the desired result.
\end{proof}

\begin{lem}
\label{lem:second_order_pole_cotransition}
Suppose that $\lambda\nearrow\mu$ in $\bY$. We have
\begin{align*}
	|\lambda|\sum_{\nu;\nu\nearrow \lambda}\frac{p^{\downarrow}(\lambda,\nu)}{(c(\mu/\lambda)-c(\lambda/\nu))^{2}}
	= \frac{1}{p^{\uparrow}(\lambda,\mu)}-1.
\end{align*}
\end{lem}
\begin{proof}
The proof is similar to that of Lemma~\ref{lem:second_order_pole_transition}.
Let us first differentiate both sides of (\ref{eq:cotrans_Cauchy_contents}) to get
\begin{align}
\label{eq:derH}
	H'_{\lambda}(z) = 1 + |\lambda|\sum_{\nu;\nu\nearrow \lambda}\frac{p^{\downarrow}(\lambda,\nu)}{(z-c(\lambda/\nu))^{2}}.
\end{align}
From (\ref{eq:CauchyInv_rational}), the derivative $H'_{\lambda}$ reads
\begin{equation*}
	H'_{\lambda}(z) = \sum_{i=1}^{d}\frac{\prod_{j=1;j\neq i}^{d}(z-\x_{i})}{\prod_{j=1}^{d-1}(z-\y_{j})}+(\text{terms vanishing at $\x_{i}$, $i=1,\dots, d$}),
\end{equation*}
where $\x_{1},\dots, \x_{d}$ and $\y_{1},\dots, \y_{d-1}$ are the local minima and maxima of the profile of $\lambda$.
As $\lambda\nearrow\mu$, we must have $c(\mu/\lambda)=\x_{i}$ for some $i$.
Evaluating the above equation at $z=\x_{i}$, we obtain
\begin{align*}
	H'_{\lambda}(c(\mu/\lambda)) = \frac{\prod_{j=1;j\neq i}^{d}(z-\x_{i})}{\prod_{j=1}^{d-1}(z-\y_{j})}=\frac{1}{p^{\uparrow}(\lambda,\mu)}.
\end{align*}
Here, the second equality follows from (\ref{eq:transition_rational}).
We conclude the desired result by evaluating (\ref{eq:derH}) at $z=c(\mu/\lambda)$.
\end{proof}

\begin{lem}
For $\lambda\nearrow\mu$, we have
\begin{equation}
\label{eq:Cauchy_transform_adding_box}
	\frac{(z-c(\mu/\lambda))^{2}}{(z-c(\mu/\lambda)-1)(z-c(\mu/\lambda)+1)}G_{\lambda}(z) = G_{\mu}(z).
\end{equation}
\end{lem}
\begin{proof}
As we can add a box to $\lambda$ to get $\mu$, $z = c(\lambda/\mu)$ is certainly a pole of $G_{\lambda}(z)$ and a zero point of $G_{\mu}(z)$.
That is why we need to multiply $(z-c(\mu/\lambda))^{2}$ in the formula.

In order to understand the division by $z-c(\mu/\lambda)\pm 1$, we first need to note that $z = c(\mu/\lambda) \pm 1$ cannot be poles of $G_{\lambda}(z)$.
Then there are four cases depending on whether each of them is a zero-point or not.
The first case is that neither of $z = c(\mu/\lambda)\pm 1$ is a zero-point of $G_{\lambda}(z)$:
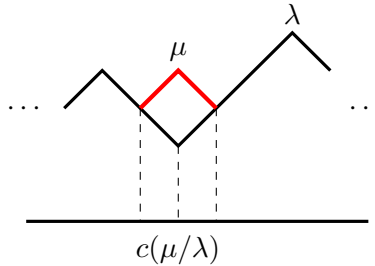
\begin{equation*}
\input{addingBox1.tikz}
\end{equation*}
In this case, they are both poles of $G_{\mu}$, thus, we have the devision by $z-c(\mu/\lambda)\pm 1$.
As another case, suppose that $z = c(\mu/\lambda)-1$ is a zero of $G_{\lambda}(z)$, but $z = c(\mu/\lambda)+1$ is not:
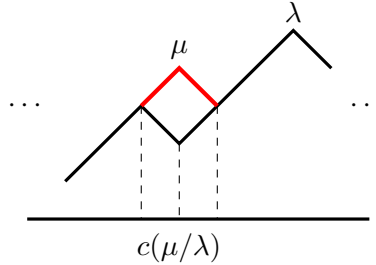
\begin{equation*}
\input{addingBox2.tikz}
\end{equation*}
Then, $z = c(\mu/\lambda)-1$ is not a zero of $G_{\mu}(z)$ and only $z = c(\mu/\lambda)+1$ is a pole,
thus the division by $z-c(\mu/\lambda)\pm 1$ occurs in this case as well.
The remaining two cases are dealt with similarly.
\end{proof}

\begin{cor}
\label{cor:transition_prob_ratio}
Suppose that $\lambda\nearrow \mu \nearrow \nu$, $\lambda\nearrow \rho\nearrow \nu$ and $\mu\neq \rho$.
Then, we have
\begin{equation*}
	\frac{(c(\nu/\mu)-c(\mu/\lambda))^{2}}{(c(\nu/\mu)-c(\mu/\lambda)-1)(c(\nu/\mu)-c(\mu/\lambda)+1)}p^{\uparrow}(\lambda,\rho) = p^{\uparrow}(\mu,\nu).
\end{equation*}
\end{cor}
\begin{proof}
Under the assumptions, we have $c(\nu/\mu)-c(\mu/\lambda)\neq \pm 1$ and $c(\nu/\mu) = c(\rho/\lambda)$.
We shall take the residue of (\ref{eq:Cauchy_transform_adding_box}) at $z=c(\nu/\mu) = c(\rho/\lambda)$ to get the desired result.
\end{proof}

\subsubsection{Shorthand notation}
\label{subsect:shorthand_contents}
In the following proofs, we will use the following notation to save the space.
When Young diagrams are labeled by superscript as $\lambda^{(i)}$ with indices $i$,
we write $c_{ij} = c(\lambda^{(i)}/\lambda^{(j)})$.
Of course, for this to make sense, we need $\lambda^{(j)}\nearrow \lambda^{(i)}$,
but that will be stated in each occasion.

\subsubsection{Proof of {\bf (R1)}}
The left-hand side is evaluated at loops of length $2$ of signature $(-,+)$, $\lambda^{(0)}\searrow \lambda^{(1)} \nearrow \lambda^{(0)}$.
Since this is the first place in this paper where we evaluate a planar tangle concretely,
we present as much detail of the computation as possible.
The value on the loop reads
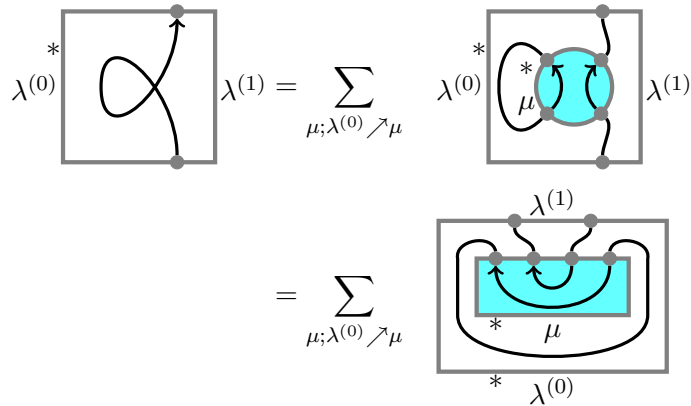
\begin{align*}
	\input{leftTurnEval.tikz}
	&= \sum_{\mu;\lambda^{(0)}\nearrow\mu}\quad \input{leftTurnEvalSum.tikz}\\
	&= \sum_{\mu;\lambda^{(0)}\nearrow\mu}\quad \input{leftTurnEvalSumIso.tikz}\quad .
\end{align*}
We isotoped the diagram into the second line because the multi-linear map is defined on planar tangles of standard shapes.
We can now start using the definition of $t_{\id}$ to obtain
\begin{align*}
	\input{leftTurnEval.tikz}
	&= \sum_{\mu;\lambda^{(0)}\nearrow\mu} \frac{1}{c(\mu/\lambda^{(0)})-c_{01}}\sqrt{\frac{f_{\Pl}(\lambda^{(1)})}{f_{\Pl}(\mu)}}\Bigg(\frac{f_{\Pl}(\mu)}{f_{\Pl}(\lambda^{(0)})}\Bigg)^{3/2} \\
	&= \sqrt{\frac{f_{\Pl}(\lambda^{(1)})}{f_{\Pl}(\lambda^{(0)})}}\sum_{\mu;\lambda^{(0)}\nearrow\mu} \frac{1}{c(\mu/\lambda^{(0)})-c_{01}}\frac{f_{\Pl}(\mu)}{f_{\Pl}(\lambda^{(0)})}.
\end{align*}
Remember the notation introduced in Section~\ref{subsect:shorthand_contents}.
Now, the sum over $\mu$ appearing above can be evaluated by using (\ref{eq:ratio_f_transition}) and (\ref{eq:Cauchy_contents}) as
\begin{align*}
	\sum_{\mu;\lambda^{(0)}\nearrow\mu} \frac{1}{c(\mu/\lambda^{(0)})-c_{01}}\frac{f_{\Pl}(\mu)}{f_{\Pl}(\lambda^{(0)})}=&\, \sum_{\mu;\lambda^{(0)}\nearrow\mu} \frac{p^{\uparrow}(\lambda^{(0)},\mu)}{c(\mu/\lambda^{(0)})-c_{01}} \\
	=&\, - G_{\lambda^{(0)}}(c_{01}).
\end{align*}
We also know that the Cauchy transform $G_{\lambda^{(0)}}$ vanishes at the contents of the removable boxes in $\lambda^{(0)}$,
in particular,
\begin{equation*}
	G_{\lambda^{(0)}}(c_{01})=0.
\end{equation*}
Therefore, the whole result is zero.

\subsubsection{Proof of {\bf (R2-1)}}
We start by the value of the left-hand side on the loop $\lambda^{(0)}\searrow \lambda^{(1)}\searrow \lambda^{(2)}\nearrow \lambda^{(3)}\nearrow \lambda^{(0)}$, which is given by the following sum over intermediate Young diagrams:
\begin{align*}
	\input{indIndCrossingsEval.tikz}= \sum_{\mu;\lambda^{(2)}\nearrow \nu\nearrow \lambda^{(0)}}\quad \input{indIndCrossingsEvalSum.tikz}\quad .
\end{align*}

We consider two separate cases.

\noindent{\bf Case I: $\lambda^{(1)} = \lambda^{(3)}$.}
There are two possibilities for $\mu$: either $\mu = \lambda^{(1)}=\lambda^{(3)}$ or the unique Young diagram $\mu^{*}$ such that $\mu^{*}\neq \lambda^{(1)}=\lambda^{(3)}$ if it exists.
Thus, the value becomes
\begin{align*}
	\input{indIndCrossingsEval.tikz}
	=&\, \input{indIndCrossingsEvalidid.tikz}\quad + \quad \input{indIndCrossingsEvalexex.tikz} \\
	=&\, \frac{1}{(c_{01}-c_{12})^{2}}\sqrt{\frac{f_{\Pl}(\lambda^{(2)})}{f_{\Pl}(\lambda^{(0)})}} 
	+\Big(1 - \frac{1}{(c_{01}-c_{12})^{2}}\Big)\sqrt{\frac{f_{\Pl}(\lambda^{(2)})}{f_{\Pl}(\lambda^{(0)})}} \\
	=&\,  \sqrt{\frac{f_{\Pl}(\lambda^{(2)})}{f_{\Pl}(\lambda^{(0)})}}
	=\, \input{indIndStraightEval.tikz}\quad .
\end{align*}
In this computation, we might seem to have assumed the existence of $\mu^{*}$,
but $t_{\ex}$ is evaluated as $0$ if there is no consistent intermediate Young diagram.
So, we end up with the same result even if $\mu^{*}$ does not exist.
We will implicitly use this type of argument, but will not repeat it explicitly in the sequel.

\noindent{\bf Case II: $\lambda^{(1)} \neq \lambda^{(3)}$.}
In this case, the intermediate Young diagram $\mu$ must coincide with either $\lambda^{(1)}$ or $\lambda^{(3)}$,
thus we get two terms
\begin{align*}
	\input{indIndCrossingsEval.tikz}
	=&\, \input{indIndCrossingsEvalidex.tikz}\quad + \quad \input{indIndCrossingsEvalexid.tikz} \\
	=&\, \frac{1}{c_{03} - c_{32}}\sqrt{1 - \frac{1}{(c_{01}-c_{12})^{2}}}\sqrt{\frac{f_{\Pl}(\lambda^{(2)})}{f_{\Pl}(\lambda^{(0)})}} \\
	&\, +\sqrt{1 - \frac{1}{(c_{01}-c_{12})^{2}}}\frac{1}{c_{01} - c_{12}}\sqrt{\frac{f_{\Pl}(\lambda^{(2)})}{f_{\Pl}(\lambda^{(0)})}}
\end{align*}
that cancel because $c_{01} = c_{32}$ and $c_{12} = c_{03}$.
Thus, we have
\begin{align*}
	\input{indIndCrossingsEval.tikz} \quad = 0 =\quad \input{indIndStraightEval.tikz}\quad.
\end{align*}

\subsubsection{Proof of {\bf (R2-2)}}
We evaluate both sides on loops $\lambda^{(0)}\searrow \lambda^{(1)}\nearrow \lambda^{(2)}\searrow \lambda^{(3)}\nearrow \lambda^{(0)}$.
From the left hand side, we get
\begin{align*}
	\input{indResCrossesEval.tikz}=
	\sum_{\mu;\lambda^{(0)},\lambda^{(2)}\nearrow \mu}\quad \input{indResCrossesEvalSum.tikz}\quad .
\end{align*}

\noindent{\bf Case I: $\lambda^{(0)}=\lambda^{(2)}$.}
In this case, we get
\begin{align*}
	\input{indResCrossesEval.tikz}
	=&\sum_{\mu;\lambda^{(0)}\nearrow \mu}\input{indResCrossesEvalidid.tikz} \\
\label{eq:indres_relation_proof_eval}
	=& \sum_{\mu;\lambda^{(0)}\nearrow \mu}\frac{1}{c(\mu/\lambda^{(0)})-c_{01}}\frac{1}{c(\mu/\lambda^{(0)})-c_{03}}
	 	\frac{f_{\Pl}(\mu)}{f_{\Pl}(\lambda^{(0)})} \\
		&\quad \cdot \sqrt{\frac{f_{\Pl}(\lambda^{(1)})}{f_{\Pl}(\lambda^{(0)})}}\sqrt{\frac{f_{\Pl}(\lambda^{(3)})}{f_{\Pl}(\lambda^{(0)})}}.
\end{align*}
Now we narrow down the following subcases

\noindent{\bf Case I-i: $\lambda^{(1)}\neq \lambda^{(3)}$.}
We can perform the partial fraction decomposition
\begin{align*}
	&\frac{1}{c(\mu/\lambda^{(0)})-c_{01}}\frac{1}{c(\mu/\lambda^{(0)})-c_{03}} \\
	=&\frac{1}{c_{01}-c_{03}}\Bigg(\frac{1}{c(\mu/\lambda^{(0)})-c_{01}} - \frac{1}{c(\mu/\lambda^{(0)})-c_{03}}\Bigg).
\end{align*}
Thus, for each $i=1,3$, the sum over $\mu$ yields
\begin{equation*}
	\sum_{\mu;\lambda^{(0)}\nearrow \mu}\frac{1}{c(\mu/\lambda^{(0)})-c_{0i}}\frac{f_{\Pl}(\mu)}{f_{\Pl}(\lambda^{(0)})}
	= \sum_{\mu;\lambda^{(0)}\nearrow \mu}\frac{p^{\uparrow}(\lambda^{(0)},\mu)}{c(\mu/\lambda^{(0)})-c_{0i}} 
	= - G_{\lambda^{(0)}}(c_{0i})
\end{equation*}
due to (\ref{eq:Cauchy_contents}), but the right-hand side vanishes as $c_{0i}$, $i=1,3$ are both zero points of $G_{\lambda^{(0)}}$.

\noindent{\bf Case I-ii: $\lambda^{(1)} = \lambda^{(3)}$.}
We may apply Lemma~\ref{lem:second_order_pole_transition} to evaluate the sum over $\mu$ and get
\begin{align*}
	\input{indResCrossesEval.tikz} = \frac{1}{|\lambda^{(0)}|p^{\downarrow}(\lambda^{(0)},\lambda^{(1)})}\frac{f_{\Pl}(\lambda^{(1)})}{f_{\Pl}(\lambda^{(0)})} = 1,
\end{align*}
where the last equality follows from (\ref{eq:ratio_f_cotransition}).

\noindent {\bf Case II: $\lambda^{(0)}\neq \lambda^{(2)}$.}
In this case, there is no consistent intermediate Young diagram $\mu$ unless $\lambda^{(1)} = \lambda^{(3)}$,
and if $\lambda^{(1)} = \lambda^{(3)}$, there exists a unique possible $\mu=\mu^{*}$.
In the case that $\lambda^{(1)}\neq\lambda^{(3)}$, the whole result is $0$, which coincides with the value of the right-hand side.
In the other case that $\lambda^{(1)}=\lambda^{(3)}$, we get
\begin{align*}
	\input{indResCrossesEval.tikz}
	=&\, \input{indResCrossesexex.tikz}\\
	=&\, \Bigl(1 - \frac{1}{(c(\nu^{*}/\lambda^{(0)})-c_{01})^{2}}\Bigl)
		\frac{f_{\Pl}(\mu^{*})}{f_{\Pl}(\lambda^{(0)})}
		\sqrt{\frac{f_{\Pl}(\lambda^{(3)})}{f_{\Pl}(\lambda^{(2)})}}
		\sqrt{\frac{f_{\Pl}(\lambda^{(1)})}{f_{\Pl}(\lambda^{(0)})}}.
\end{align*}
We may apply Corollary~\ref{cor:transition_prob_ratio} to get
\begin{align*}
	\Bigl(1 - \frac{1}{(c(\nu^{*}/\lambda^{(0)})-c_{01})^{2}}\Bigl)
		\frac{f_{\Pl}(\mu^{*})}{f_{\Pl}(\lambda^{(0)})}
	=&\, \Bigl(1 - \frac{1}{(c(\nu^{*}/\lambda^{(0)})-c_{01})^{2}}\Bigl)p^{\uparrow}(\lambda^{(0)},\mu) \\
	=&\, p^{\uparrow}(\lambda^{(1)},\lambda^{(2)}) \\
	=&\, \frac{f_{\Pl}(\lambda^{(2)})}{f_{\Pl}(\lambda^{(1)})}.
\end{align*}
Thus, we have
\begin{align*}
	\input{indResCrossesEval.tikz}
	=\frac{f_{\Pl}(\lambda^{(2)})}{f_{\Pl}(\lambda^{(1)})}
		\sqrt{\frac{f_{\Pl}(\lambda^{(3)})}{f_{\Pl}(\lambda^{(2)})}}
		\sqrt{\frac{f_{\Pl}(\lambda^{(1)})}{f_{\Pl}(\lambda^{(0)})}}
	=\sqrt{\frac{f_{\Pl}(\lambda^{(2)})}{f_{\Pl}(\lambda^{(0)})}}
\end{align*}
as desired.

\subsubsection{Proof of {\bf (R2-3)}}
The line of proof is very similar to the previous one {\bf (R2-2)}.
As always, let us start evaluating the left-hand side on the loop
$\lambda^{(0)}\nearrow \lambda^{(1)}\searrow \lambda^{(2)}\nearrow \lambda^{(3)}\searrow \lambda^{(0)}$
and get
\begin{align*}
	\input{resIndCrossesEval.tikz}=\sum_{\mu; \mu\nearrow \lambda^{(0)}, \lambda^{(2)}}\quad \input{resIndCrossesEvalSum.tikz}\quad .
\end{align*}

\noindent{\bf Case I: $\lambda^{(0)}=\lambda^{(2)}$.}
In this case, the nontrivial sum over the intermediate Young diagrams $\mu$ remains
\begin{align*}
	\input{resIndCrossesEval.tikz}
	=&\, \sum_{\mu; \mu\nearrow \lambda^{(0)}}\quad \input{resIndCrossesidid.tikz}\\
	=&\, \sum_{\mu; \mu\nearrow \lambda^{(0)}}\frac{1}{c_{10}-c(\lambda^{(0)}/\mu)}\frac{1}{c_{30}-c(\lambda^{(0)}/\mu)}
		\frac{f_{\Pl}(\mu)}{f_{\Pl}(\lambda^{(0)})} \\
	&\quad\cdot \sqrt{\frac{f_{\Pl}(\lambda^{(1)})}{f_{\Pl}(\lambda^{(0)})}}
				\sqrt{\frac{f_{\Pl}(\lambda^{(3)})}{f_{\Pl}(\lambda^{(0)})}}
\end{align*}

We consider the following two subcases.

\noindent{\bf Case I-i: $\lambda^{(1)}\neq\lambda^{(3)}$.}
We have the partial fraction decomposition
\begin{align*}
	&\,\frac{1}{c_{10}-c(\lambda^{(0)}/\mu)}\frac{1}{c_{30}-c(\lambda^{(0)}/\mu)} \\
	=&\, \frac{1}{c_{30}-c_{10}}
	\Bigl(\frac{1}{c_{10}-c(\lambda^{(0)}/\mu)} - \frac{1}{c_{30}-c(\lambda^{(0)}/\mu)}\Bigr).
\end{align*}
From (\ref{eq:cotrans_Cauchy_contents}), we can see that, for each $i=0,3$,
\begin{align*}
	\sum_{\mu; \mu\nearrow \lambda^{(0)}}\frac{1}{c_{i0}-c(\lambda^{(0)}/\mu)}\frac{f_{\Pl}(\mu)}{f_{\Pl}(\lambda^{(0)})} 
	= \sum_{\mu; \mu\nearrow \lambda^{(0)}}\frac{|\lambda^{(0)}|p^{\downarrow}(\lambda^{(0)},\mu)}{c_{i0}-c(\lambda^{(0)}/\mu)}
	= c_{i0}
\end{align*}
as $c_{i0}$, $i=1,3$ are zero-points of $H_{\lambda^{(0)}}$.
Therefore, the object of interest comes down to
\begin{align*}
	\input{resIndCrossesEval.tikz}
	=&\, \frac{1}{c_{30}-c_{10}}(c_{10}-c_{30})\sqrt{\frac{f_{\Pl}(\lambda^{(1)})}{f_{\Pl}(\lambda^{(0)})}}
				\sqrt{\frac{f_{\Pl}(\lambda^{(3)})}{f_{\Pl}(\lambda^{(0)})}} \\
	=&\, -\sqrt{\frac{f_{\Pl}(\lambda^{(1)})}{f_{\Pl}(\lambda^{(0)})}}\sqrt{\frac{f_{\Pl}(\lambda^{(3)})}{f_{\Pl}(\lambda^{(0)})}},
\end{align*}
which is certainly the desired result in this case.

\noindent{\bf Case I-ii: $\lambda^{(1)}=\lambda^{(3)}$.}
In this case, we use Lemma~\ref{lem:second_order_pole_cotransition} to deal with the sum over intermediate Young diagrams:
\begin{align*}
	\input{resIndCrossesEval.tikz}
	=&\, |\lambda^{(0)}|\sum_{\mu; \mu\nearrow \lambda^{(0)}}\frac{p^{\downarrow}(\lambda^{(0)},\mu)}{(c_{10}-c(\lambda^{(0)}/\mu))^{2}} \frac{f_{\Pl}(\lambda^{(1)})}{f_{\Pl}(\lambda^{(0)})} \\
	=&\, \Bigg(\frac{1}{p^{\uparrow}(\lambda^{(0)},\lambda^{(1)})}-1\Bigg)\frac{f_{\Pl}(\lambda^{(1)})}{f_{\Pl}(\lambda^{(0)})} \\
	=&\, 1- \frac{f_{\Pl}(\lambda^{(1)})}{f_{\Pl}(\lambda^{(0)})}.
\end{align*}
This coincides with the corresponding value of the right-hand side.

\noindent{\bf Case II: $\lambda^{(0)}\neq\lambda^{(2)}$.}
In this case, there is no consistent $\mu$ unless $\lambda^{(1)}=\lambda^{(3)}$,
thus the result is $0$ if $\lambda^{(1)}\neq \lambda^{(3)}$, which coincides with the value of the right-hand side.
If $\lambda^{(1)}=\lambda^{(3)}$, there is a unique consistent Young diagram $\mu=\mu^{*}$ that gives
\begin{align*}
	\input{resIndCrossesEval.tikz}
	=&\, \input{resIndCrossesexex.tikz} \\
	=&\, \Big(1-\frac{1}{(c_{10}-c(\lambda^{(0)}/\mu^{*}))^{2}}\Big)\frac{f_{\Pl}(\lambda^{(1)})}{f_{\Pl}(\lambda^{(0)})}\frac{f_{\Pl}(\mu^{*})}{f_{\Pl}(\lambda^{(2)})}\sqrt{\frac{f_{\Pl}(\lambda^{(2)})}{f_{\Pl}(\lambda^{(0)})}}
\end{align*}
Similarly to {\bf Case II} of the previous {\bf (R2-2)}, we obtain from Corollary~\ref{cor:transition_prob_ratio} 
\begin{align*}
	\Big(1-\frac{1}{(c_{10}-c(\lambda^{(0)}/\mu^{*}))^{2}}\Big)\frac{f_{\Pl}(\lambda^{(1)})}{f_{\Pl}(\lambda^{(0)})}=\frac{f_{\Pl}(\lambda^{(2)})}{f_{\Pl}(\mu^{*})}
\end{align*}
that completes the proof.

\subsubsection{Proof of {\bf (R3)}}
Both sides are evaluated on loops of the form $\lambda^{(0)}\searrow \lambda^{(1)}\searrow \lambda^{(2)}\nearrow \lambda^{(3)}\nearrow \lambda^{(4)}\nearrow \lambda^{(5)}\nearrow \lambda^{(0)}$.
The goal is to show
\begin{align}
\label{eq:YBEStateSum}
	\input{yangBaxterLeftEval.tikz}\quad =\quad \input{yangBaxterRightEval.tikz}\quad .
\end{align}

\noindent{\bf Case I: $\lambda^{(1)} = \lambda^{(5)}$ and $\lambda^{(2)} = \lambda^{(4)}$.} 
This is the case when the paths $\lambda^{(3)}\nearrow \lambda^{(2)} \nearrow \lambda^{(1)}\nearrow \lambda^{(0)}$ and $\lambda^{(3)}\nearrow \lambda^{(4)} \nearrow \lambda^{(5)}\nearrow \lambda^{(0)}$ are identical.
The left-hand side of (\ref{eq:YBEStateSum}) yields
\begin{align*}
	&\, \input{yangBaxterLeftididid.tikz} + \input{yangBaxterLeftexidex.tikz} \\
	=&\, \Biggl( \frac{1}{c_{01}-c_{12}}\frac{1}{(c_{12} - c_{23})^{2}}
	+\frac{1}{c_{01}-c_{23}}\Bigl(1 - \frac{1}{(c_{12} - c_{23})^{2}}\Bigr)\Biggr) \sqrt{\frac{f_{\Pl}(\lambda^{(3)})}{f_{\Pl}(\lambda^{(0)})}},	
\end{align*}
whereas the right-hand side of (\ref{eq:YBEStateSum}) does
\begin{align*}
	&\, \input{yangBaxterRightididid.tikz} + \input{yangBaxterRightexidex.tikz} \\
	=&\, \Biggl( \frac{1}{(c_{01} - c_{12})^{2}}\frac{1}{c_{12} - c_{23}}
	+\Bigl(1-\frac{1}{(c_{01}-c_{12})^{2}}\Bigr)\frac{1}{c_{01}-c_{23}}\Biggr)\sqrt{\frac{f_{\Pl}(\lambda^{(3)})}{f_{\Pl}(\lambda^{(0)})}}.
\end{align*}
It can be readily checked that these two expressions coincide.

\noindent{\bf Case II: $\lambda^{(1)} = \lambda^{(5)}$ and $\lambda^{(2)}\neq\lambda^{(4)}$.}
In this case, we get from the left-hand side of (\ref{eq:YBEStateSum})
\begin{align*}
	&\, \input{yangBaxterLeftexidid.tikz}+\input{yangBaxterLeftididex.tikz} \\
	=& \Bigg(\frac{1}{c_{01}-c_{12}}\frac{1}{c_{12}-c_{23}}
	+\frac{1}{c_{01}-c_{23}}\frac{1}{c_{23}-c_{12}}\Biggr)\sqrt{1-\frac{1}{(c_{12}-c_{23})^{2}}} \sqrt{\frac{f_{\Pl}(\lambda^{(3)})}{f_{\Pl}(\lambda^{(0)})}},
\end{align*}
and from the right-hand side
\begin{align*}
	\input{yangBaxterRightidexid.tikz} = \frac{1}{c_{01}-c_{12}}\frac{1}{c_{01}-c_{23}}\sqrt{1-\frac{1}{(c_{12}-c_{23})^{2}}}\sqrt{\frac{f_{\Pl}(\lambda^{(3)})}{f_{\Pl}(\lambda^{(0)})}}.
\end{align*}
The coincidence of the two values can be checked.

\noindent{\bf Case III: $\lambda^{(1)}\neq \lambda^{(5)}$ and $\lambda^{(2)} = \lambda^{(4)}$.}
We omit this case as it exhibits a similar pattern as {\bf Case II}.

\noindent{\bf Case IV: $\lambda^{(1)}\neq \lambda^{(5)}$ and $\lambda^{(2)}\neq \lambda^{(4)}$.}
This case breaks down to the following three subcases:
\begin{itemize}
\item {\bf Case IV-i: $\lambda^{(0)}/\lambda^{(1)}=\lambda^{(5)}/\lambda^{(4)}$, $\lambda^{(1)}/\lambda^{(2)}=\lambda^{(4)}/\lambda^{(3)}$, and $\lambda^{(2)}/\lambda^{(3)}=\lambda^{(0)}/\lambda^{(5)}$.}
\item {\bf Case IV-ii: $\lambda^{(0)}/\lambda^{(1)}=\lambda^{(4)}/\lambda^{(3)}$, $\lambda^{(1)}/\lambda^{(2)}=\lambda^{(0)}/\lambda^{(5)}$, and $\lambda^{(2)}/\lambda^{(3)}=\lambda^{(5)}/\lambda^{(4)}$.}
\item {\bf Case IV-iii: $\lambda^{(0)}/\lambda^{(1)}=\lambda^{(4)}/\lambda^{(3)}$, $\lambda^{(1)}/\lambda^{(2)}=\lambda^{(5)}/\lambda^{(4)}$, and $\lambda^{(2)}/\lambda^{(3)}=\lambda^{(0)}/\lambda^{(5)}$.}
\end{itemize}
All these three cases are dealt with similarly; for each case, the left-hand and right-hand sides of (\ref{eq:YBEStateSum}) have only one term.
For example in {\bf Case IV-i}, we have
\begin{align*}
	\input{yangBaxterLeftidexex.tikz}=\input{yangBaxterRightexexid.tikz}
\end{align*}
since both sides give
\begin{align*}
	\sqrt{1-\frac{1}{(c_{01}-c_{23})^{2}}}\sqrt{1-\frac{1}{(c_{12}-c_{23})^{2}}}\frac{1}{c_{01}-c_{12}}.
\end{align*}

\subsection{Moments and Boolean cumulants}
Recall that elements in $P_{\ul{\emptyset}}$ are considered function on $\bY$.
Among others, we shall see those that are identified with moments and Boolean cumulants,
which we introduced in Section~\ref{sect:def_moments_cumulants}.

We first introduce the right-turn element depicted by a black dot as follows:
\begin{align*}
	\input{JME.tikz}\, :=\, \input{rightTurn.tikz}\, \in P_{(-,+)}.
\end{align*}
Then, for $k\in\bZ_{\geq 0}$, the $k$-th {\it power} of the right-turn element is denoted by
\begin{align*}
	\input{JMEpower.tikz}\, :=\, \input{JMEpile.tikz}\, ,
\end{align*}
where, in the right-hand side, we have omitted the external disks for the black dots,
but this should not be confusing as long as we remember the location of $*$ relative to each black dot.
The $k=0$ case is understood as a single straight line.

\begin{prop}
\label{prop:JME_value}
For $\mu\nearrow\lambda$, we have
\begin{align*}
	\input{JMEEval.tikz}\, = c(\lambda/\mu)\sqrt{\frac{f_{\Pl}(\mu)}{f_{\Pl}(\lambda)}}.
\end{align*}
\end{prop}
\begin{proof}
The result follows from the following direct calculation:
\begin{align*}
	\input{JMEEval.tikz}\, 
	=&\, \sum_{\nu;\nu\nearrow\mu}\input{rightTurnEvalShade.tikz}\\
	=&\, \sum_{\nu;\nu\nearrow\mu}\frac{1}{c(\lambda/\mu)-c(\mu/\nu)}\frac{f_{\Pl}(\nu)}{f_{\Pl}(\mu)}\sqrt{\frac{f_{\Pl}(\mu)}{f_{\Pl}(\lambda)}} \\
	=&\, |\mu|\sum_{\nu;\nu\nearrow\mu}\frac{p^{\downarrow}(\mu,\nu)}{c(\lambda/\mu)-c(\mu/\nu)}\sqrt{\frac{f_{\Pl}(\mu)}{f_{\Pl}(\lambda)}} \\
	=&\, c(\lambda/\mu)\sqrt{\frac{f_{\Pl}(\mu)}{f_{\Pl}(\lambda)}}.
\end{align*}
Here, we used (\ref{eq:cotrans_Cauchy_contents}) in the last equality and the fact that $c(\lambda/\mu)$ is a zero-point of $H_{\mu}$.
\end{proof}

In the following theorems, we omit the external disk for elements in $P_{(\emptyset)}$ (recall Section~\ref{sect:planar_alg_convention}).
Let us start looking at the moments of Young diagrams.

\begin{thm}
\label{thm:moment_diagram}
For $k\in\bZ_{>0}$, we have
\begin{align*}
	M_{k} = \input{moment.tikz}.
\end{align*}
\end{thm}
\begin{proof}
We first compute the value of the right-hand side at $\lambda\in \bY$:
\begin{align*}
	\input{momentEval.tikz}\, &=\sum_{\mu;\lambda\nearrow \mu}\, \input{momentEvalInt.tikz}.
\end{align*}
Here, we can use Proposition~\ref{prop:JME_value} to evaluate the right-hand side and get
\begin{align*}
	\input{momentEval.tikz}\, &=\sum_{\mu;\lambda\nearrow \mu}c(\mu/\lambda)^{k}\, \frac{f_{\Pl}(\mu)}{f_{\Pl}(\lambda)}\\
	&=\sum_{\mu;\lambda\nearrow \mu}c(\mu/\lambda)^{k}\, p^{\uparrow}(\lambda,\mu).
\end{align*}
Now, we expand (\ref{eq:Cauchy_contents}) in terms of a large enough $z$:
\begin{align*}
	G_{\lambda}(z) &= z + \sum_{k=1}^{\infty}z^{-k-1}\sum_{\mu;\lambda\nearrow\mu}c(\mu/\lambda)^{k}\, p^{\uparrow}(\lambda,\mu) \\
	&=z + \sum_{k=1}^{\infty}z^{-k-1}\, \input{momentEval.tikz}.
\end{align*}
When we compare this expansion with (\ref{eq:def_moments}), we see the desired result.
\end{proof}

The Boolean cumulants are realized by the following elements in $P_{\ul{\emptyset}}$.
\begin{thm}
For $k\in \bZ_{\geq 0}$, we have
\begin{align*}
	B_{k+2} = \input{boolCum.tikz}\quad .
\end{align*}
\end{thm}
\begin{proof}
The proof is very similar to that of Theorem~\ref{thm:moment_diagram}.
Let us evaluate the right-hand side on $\lambda\in \bY$:
\begin{align*}
	\input{boolCumEval.tikz}\, &=\sum_{\mu;\mu\nearrow\lambda}\, \input{boolCumEvalInt.tikz} \\
	&=\sum_{\mu;\mu\nearrow\lambda}c(\lambda/\mu)^{k}\, \frac{f_{\Pl}(\mu)}{f_{\Pl}(\lambda)} \\
	&=|\lambda|\sum_{\mu;\mu\nearrow\lambda}c(\lambda/\mu)^{k}\, p^{\downarrow}(\lambda,\mu),
\end{align*}
where we used Proposition~\ref{prop:JME_value} to go to the second line.
When we expand (\ref{eq:cotrans_Cauchy_contents}) in a sufficiently large $z$, we have
\begin{align*}
	H_{\lambda}(z) &= z -\sum_{k=0}^{\infty}z^{-k-1}|\lambda|\sum_{\mu;\mu\nearrow\lambda}c(\lambda/\mu)^{k}\, p^{\downarrow}(\lambda,\mu) \\
	&= z - \sum_{k=0}^{\infty}z^{-k-1}\quad \input{boolCumEval.tikz}
\end{align*}
Comparing this with the definition of Boolean cumulants (\ref{eq:def_Boolean_cumulants}), we get the desired result.
\end{proof}

\subsection{Normalized characters}
\label{sect:normalizedCharacters}
For $k\in \bZ_{>0}$, we introduce the following element in $P_{(-^{k},+^{k})}$, where $(-^{k},+^{k}) = (\underbrace{-,\dots,-}_{k},\underbrace{+,\dots,+}_{k})$:
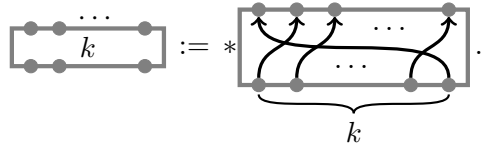
\begin{align*}
	\input{cycleBox.tikz}\, := \, \input{cycleString.tikz}\, .
\end{align*}
While, in the left-hand side, the distinguished boundary component is implicit,
the notation should not be confusing in practice as long as we understand it in $P_{(-^{k},+^{k})}$.

%We use the following shorthand notation: $(-^{k},+^{k}) = (\underbrace{-,\dots,-}_{k},\underbrace{+,\dots,+}_{k})$ with $k\in\bZ_{>0}$.

\begin{thm}
\label{thm:characterDiagram}
Let $\pi=(\pi_{1},\dots,\pi_{l})$ be a partition. Then, we have
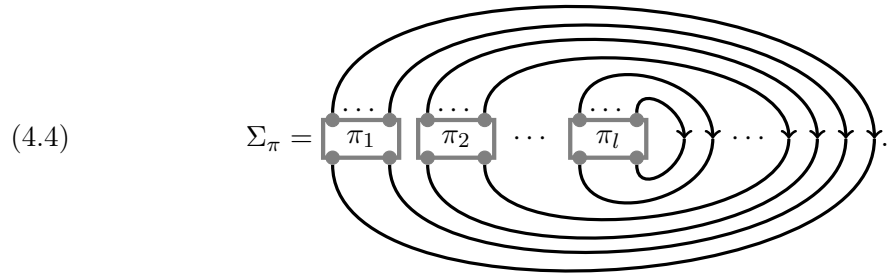
\begin{align}
\label{eq:diagramForCharacter}
	\Sigma_{\pi} = \input{characterDiagram.tikz}\, .
\end{align}
\end{thm}

To complete the proof, we will need some representation theory of symmetric groups, but we will not dig into it in this section.
Instead, we present our proof only by means of the planar algebra until the point we need representation theory,
and postpone further details to Appendix~\ref{app:reptheorSG}.

\begin{proof}
Let $k\in\bZ_{>0}$ and $h\in P_{\ul{\emptyset}}$.
We first study the function
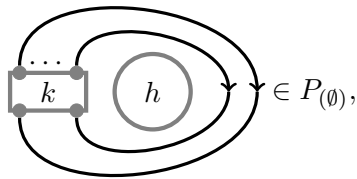
\begin{equation*}
	\input{cycleWithGeneric.tikz}\, \in P_{(\emptyset)},
\end{equation*}
where, as usual, we omit the external disk.
Placing $h$ inside the internal disk means that $h$ is evaluated there.

Let us evaluate the above function at $\lambda\in\bY_{n}$ with $n\in\bN$:
\begin{align*}
	&\,\input{cycleWithGenericEval.tikz}\\
	=&\,\sum_{\substack{\lambda^{(1)},\dots,\lambda^{(k)}\\ \lambda^{(k)}\nearrow\cdots\nearrow \lambda^{(1)}\nearrow\lambda}}\,
	\input{cycleWithGenericEvalStates.tikz} \\
	=&\, \sum_{\substack{\lambda^{(1)},\dots,\lambda^{(k)}\\ \lambda^{(k)}\nearrow\cdots\nearrow \lambda^{(1)}\nearrow\lambda}}\prod_{i=1}^{k-1}\frac{1}{c(\lambda^{(i-1)}/\lambda^{(i)})-c(\lambda^{(i)}/\lambda^{(i+1)})}\frac{f_{\Pl}(\lambda^{(k)})}{f_{\Pl}(\lambda)}h(\lambda^{(k)}),
\end{align*}
where we agree that $\lambda^{(0)}=\lambda$.
Applying this computation recursively allows us to evaluate the right-hand side of (\ref{eq:diagramForCharacter}) and get
\begin{align}
\label{eq:characterDiagramEval}
	&\, \input{characterDiagramEval.tikz} \\
	&\,=\sum_{\substack{\lambda^{(1)},\dots,\lambda^{(k)}\\ \lambda^{(k)}\nearrow\cdots\nearrow \lambda^{(1)}\nearrow\lambda}}\prod_{i=1}^{l}\prod_{j=\pi_{1}+\cdots+\pi_{i-1}+1}^{\pi_{1}+\cdots+\pi_{i}-1}\frac{1}{c(\lambda^{(j-1)}/\lambda^{(j)})-c(\lambda^{(j)}/\lambda^{(j+1)})}\frac{f_{\Pl}(\lambda^{(k)})}{f_{\Pl}(\lambda)}. \notag
\end{align}
From (\ref{eq:ratio_f_cotransition}) and (\ref{eq:cotrans_dim_ratio}), we can see that
\begin{equation*}
	\frac{f_{\Pl}(\lambda^{(k)})}{f_{\Pl}(\lambda^{(0)})} = (n \downharpoonright k)\frac{\dim\lambda^{(k)}}{\dim\lambda}.
\end{equation*}
We will see (\ref{eq:character_from_representation}) in Appendix~\ref{app:reptheorSG} that turns the right-hand side of (\ref{eq:characterDiagramEval}) to
\begin{align*}
	(n\downharpoonright k)\frac{\chi^{\lambda}_{\pi\cup (1^{n-k})}}{\dim \lambda}
\end{align*}
when $k\leq n$.
When $k>n$, the sum in (\ref{eq:characterDiagramEval}) is empty yielding $0$.
Comparing this with the definition of normalized characters (\ref{eq:defNormalizedCharacter}),
we obtain the desired result.
\end{proof}

\appendix

\section{Representations of symmetric groups}
\label{app:reptheorSG}
In this appendix, we briefly recall the representation theory of symmetric groups.
We refer to~\cite{Sagan2001} for a general background for symmetric groups and their representations.
For our purpose, the picture given in~\cite{VershikOkounkov2005} is particularly useful.
Although one of the upshots in~\cite{VershikOkounkov2005} is that the irreducible representations of a symmetric group $S_{n}$ are labelled by $\bY_{n}$, we take that part for granted and start discussing branching right away.

\subsection{Branching of representations}
We inherit the notation from Section~\ref{subsect:normalizedCharacter};
for each $\lambda\in\bY$, we write $V^{\lambda}$ for the representation space of an irreducible representation labelled by $\lambda$.
When $\lambda\in \bY_{n}$, $V^{\lambda}$ is an irreducible representation of $S_{n}$, but decomposes as a representation of $S_{n-1}$ as follows:
\begin{equation*}
    V^{\lambda} = \bigoplus_{\mu\in\bY_{n-1};\mu\nearrow\lambda}V^{\mu}.
\end{equation*}
We may continue this procedure until $S_{0}$ to get a decomposition into one-dimensional vector spaces:
\begin{equation*}
    V^{\lambda} = \bigoplus_{\varpi:\emptyset\nearrow \lambda^{(1)}\nearrow\cdots\nearrow \lambda^{(n)}=\lambda} V_{\varpi},
\end{equation*}
where $V_{\varpi}$ are copies of $V^{\emptyset}$ labelled by paths $\varpi$ in the Young graph from $\emptyset$ to $\lambda$.
Each of such paths are identified with a standard Young tableau of shape $\lambda$. We write $\ST (\lambda)$ for the set of standard tableaux of shape $\lambda$.
When we fix a non-zero vector $v_{\varpi}\in V_{\varpi}$, the vectors $\{v_{\varpi}\}_{\varpi\in\ST(\lambda)}$ form a basis of $V^{\lambda}$, called the Gelfand--Tsetlin (GZ) basis.

The GZ basis is not unique because we can choose non-zero scalar multiples of the basis vectors.
It has been shown in~\cite{VershikOkounkov2005} that there is a uniform choice of basis vectors so that
the matrix element formulae in the next section hold.

\subsection{Matrix elements of simple transpose}
\label{sect:app_matrix_elements}
Let $t_{i}=(i\, i+1)$ be a simple transpose. We can write down an explicit formula for the matrix elements of each $t_{i}$ in the GZ basis as follows.

Suppose that $\varpi\in \ST(\lambda)$ with $\lambda\in \bY_{n}$ is given by
\begin{align*}
	\varpi : \emptyset = \lambda^{(0)}\nearrow \cdots \nearrow \lambda^{(n)}= \lambda.
\end{align*}
For each $i=1,\dots, n-1$, let us write down the vector $t_{i}v_{\varpi}$.
There are two possible cases.
The first case is that $c(\lambda^{(i+1)}/\lambda^{(i)})-c(\lambda^{(i)}/\lambda^{(i-1)})=\pm 1$.
In this case, $t_{i}v_{\varpi}$ is proportional to $v_{\varpi}$ and is given by
\begin{align*}
	t_{i}v_{\varpi} = \pm v_{\varpi} = \frac{1}{c(\lambda^{(i+1)}/\lambda^{(i)})-c(\lambda^{(i)}/\lambda^{(i-1)})}v_{\varpi}.
\end{align*}
Writing the right-hand side in this way seems redundant as the coefficient is just $\pm 1$, but it allows for a compact formula combined with the second case.

The second case is that $c(\lambda^{(i+1)}/\lambda^{(i)})-c(\lambda^{(i)}/\lambda^{(i-1)})\neq \pm 1$.
In this case, there is an another path
\begin{align*}
	\varpi' : \emptyset = \lambda^{(0)}\nearrow \cdots \nearrow \lambda^{(i-1)} \nearrow \lambda^{(i)'} \nearrow \lambda^{(i+1)}\nearrow \cdots \nearrow \lambda^{(n)}=\lambda
\end{align*}
such that
\begin{align*}
	\lambda^{(i+1)}/\lambda^{(i)} = \lambda^{(i)'}/\lambda^{(i-1)},\quad \lambda^{(i)}/\lambda^{(i-1)} = \lambda^{(i+1)}/\lambda^{(i)'}.
\end{align*}
The two dimensional space $\bC v_{\varpi}\oplus \bC v_{\varpi'}$ is invariant under $t_{i}$, and in that basis, $t_{i}$ takes the matrix form
\begin{align*}
	t_{i} = \left(
	\begin{array}{cc}
	r^{-1} & \sqrt{1-r^{-2}} \\
	\sqrt{1-r^{-2}} & -r^{-1}
	\end{array}
	\right),\quad r = c(\lambda^{(i+1)}/\lambda^{(i)})-c(\lambda^{(i)}/\lambda^{(i-1)})
\end{align*}
In either case, the diagonal matrix element of $t_{i}$ is always of the form
\begin{align*}
	\frac{1}{c(\lambda^{(i+1)}/\lambda^{(i)})-c(\lambda^{(i)}/\lambda^{(i-1)})}
\end{align*}
when acted on $v_{\varpi}$.

\subsection{Characters}
For $\lambda\in \bY_{n}$ and $\pi\in \cP_{k}$ with $k\leq n$, the character $\chi^{\lambda}_{\pi\cup (1^{n-k})}$ is computed as
\begin{align*}
	\chi^{\lambda}_{\pi\cup (1^{n-k})} = \Tr_{V^{\lambda}}(\sigma),
\end{align*}
where $\sigma\in S_{n}$ has the cycle type $\pi\cup (1^{n-k})$.
We can choose such $\sigma$ as
\begin{align*}
	\sigma =
	&(n\quad n-1\quad \dots\quad n-\pi_{1}+1) \\
	&(n-\pi_{1}\quad n-\pi_{1}-1\quad \dots\quad n-\pi_{1}-\pi_{2}+1) \\
	&\quad \vdots \\
	&(n-k+\pi_{l} \quad n-k+\pi_{k}-1\quad\dots\quad n-k+1)
\end{align*}
so that it acts on the last $k$ letters.

Recall that $V^{\lambda}$ has the GZ-basis $\{v_{\varpi}\}_{\varpi\in\ST(\lambda)}$.
For each
\begin{equation*}
	\varpi : \emptyset = \lambda^{(n)}\nearrow \cdots \nearrow \lambda^{(1)}\nearrow \lambda^{(0)}=\lambda,
\end{equation*}
the diagonal matrix element of $\sigma$ on $v_{\varpi}$, denoted by $M_{\varpi}[\sigma]$, can be computed as
\begin{align*}
	M_{\varpi}[\sigma] = \prod_{i=1}^{l}\prod_{j=\pi_{1}+\cdots+\pi_{i-1}+1}^{\pi_{1}+\cdots+\pi_{i}-1}\frac{1}{c(\lambda^{(j-1)}/\lambda^{(j)})-c(\lambda^{(j)}/\lambda^{(j+1)})}
\end{align*}
using the results in Section~\ref{sect:app_matrix_elements}.
Note that we have changed the labelling of the Young diagrams along a path to the decreasing way
whereas they were labelled in the increasing way in Section~\ref{sect:app_matrix_elements}.
The matrix element $M_{\varpi}[\sigma]$ is obviously independent of the first $n-k$ part
\begin{align*}
	\emptyset = \lambda^{(n)}\nearrow \cdots \nearrow \lambda^{(k+1)}
\end{align*}
of the path.
In other words, fixing the last $k$ part
\begin{align*}
	\lambda^{(k)}\nearrow \cdots \nearrow \lambda^{(0)}=\lambda,
\end{align*}
there are $\# \ST(\lambda^{(k)})$ number of paths $\varpi$ that give the same diagonal matrix element $M_{\varpi}[\sigma]$.
Therefore, we can compute the character as
\begin{align}
\label{eq:character_from_representation}
	\chi^{\lambda}_{\pi\cup (1^{n-k})}
	 =\sum_{\substack{\lambda^{(1)},\dots,\lambda^{(k)}\\ \lambda^{(k)}\nearrow\cdots\nearrow\lambda^{(1)}\nearrow\lambda}}\dim\lambda^{(k)}\, \prod_{i=1}^{l}\prod_{j=\pi_{1}+\cdots+\pi_{i-1}+1}^{\pi_{1}+\cdots+\pi_{i}-1}\frac{1}{c(\lambda^{(j-1)}/\lambda^{(j)})-c(\lambda^{(j)}/\lambda^{(j+1)})},
\end{align}
where we substituted $\#\ST(\lambda^{(k)}) = \dim \lambda^{(k)}$.

\section{Towards an alternative proof of Theorem~\ref{thm:characterDiagram}}
\label{app:charactersAlter}
In Theorem~\ref{thm:characterDiagram}, we have identified certain diagrams with the normalized characters,
where the proof went by using the definition of the planar algebra and representation theory of symmetric groups.
Although we have explicitly written down the trace of a symmetric group action
to get a formula for a character in (\ref{eq:character_from_representation}), there is another common way to evaluate it.

The Schur--Frobenius formula reads
\begin{align}
\label{eq:SchurFrobenius}
	p_{\pi} = \sum_{\lambda\in\bY_{n}}\chi^{\lambda}_{\pi}s_{\lambda},\quad \pi\in \cP_{n},\, n\in\bN,
\end{align}
where $s_{\lambda}$ and $p_{\pi}$ are Schur and power-sum symmetric functions (see \cite{Fulton1996,macdonald1998symmetric})
and the characters appear as expansion coefficients.
This formula connects the representation theory of symmetric groups to the theory of symmetric functions,
and allows us to compute the characters in terms of symmetric functions.
In particular, the Frobenius character formula provides, for $k\in \bN$,
\begin{align}
\label{eq:FrobeniusCharFormula}
	\Sigma_{(k)}(\lambda) = -\frac{1}{k}\int_{\Gamma} H_{\lambda}(z)H_{\lambda}(z-1)\cdots H_{\lambda}(z-k+1)dz,
\end{align}
where $\Gamma$ is an integral contour in $\bC$ that encloses all the poles of the integrand,
and we have omit the factor $\frac{1}{2\pi\bbmi}$ that always comes with a complex integral.
The formula (\ref{eq:FrobeniusCharFormula}) is the version found in~\cite{Biane2003}.
It only covers the case of single part partitions, but was extended in~\cite{RS2008} to general partitions.

On the contrary, we can even use (\ref{eq:SchurFrobenius}) as a {\it definition} of characters,
and as such, the notion of characters becomes independent of representation theory but
intrinsic in the theory of symmetric functions.
In fact, this is how one defines the Jack deformed characters that have so far no representation theoretical origin except at special parameters.

As we have discussed in Introduction, planar algebras for the Young graph seem to admit a Jack deformation.
It would be therefore reasonable to expect that everything works without representation theory even in the undeformed case
considered in this paper.

The purpose of this Appendix is to see a possibility of proving Theorem~\ref{thm:characterDiagram} directly comparing the diagrams with (\ref{eq:FrobeniusCharFormula}),
In fact, the proof of Theorem~\ref{thm:characterDiagram} was the only place where we needed representation theory.

\subsection{Frobenius character formula, revisited}
Let us first observe another form of (\ref{eq:FrobeniusCharFormula}).
For $\lambda\in \bY$ and $n\in \bZ_{>0}$, we set a function
\begin{equation}
\label{eq:defFunctionF}
	F_{\lambda}^{(n)} (z_{1},\dots, z_{n}) = \prod_{i=1}^{n-1}\frac{1}{z_{i}-z_{i+1}}\prod_{1\leq i<j\leq n}\frac{(z_{i}-z_{j})^{2}}{(z_{i}-z_{j}-1)(z_{i}-z_{j}+1)}\prod_{i=1}^{n}H_{\lambda}(z_{i})
\end{equation}
of $z_{1},\dots, z_{n}$.
Here, we understand $F_{\lambda}^{(1)}(z) = H_{\lambda}(z)$.
This function has poles in each variable at the local maxima of the profile of $\lambda$
as well as at the points $z_{i}-z_{j}=\pm 1$ for $i\neq j$.

\begin{thm}
\label{thm:FrobeniusSatelite}
Let $\lambda\in\bY$ and $n\in\bN$.
We define $I_{\lambda}^{(k)}(z_{k+1},\dots,z_{n})$, $k=0,\dots, n-1$ recursively as follows:
\begin{align*}
	I_{\lambda}^{(0)}(z_{1},\dots,z_{n})&= F_{\lambda}^{(n)}(z_{1},\dots, z_{n}), \\
	I_{\lambda}^{(k)}(z_{k+1},\dots, z_{n})&= \int_{\Gamma_{z_{k+1}}^{(k)}}I_{\lambda}^{(k-1)}(z_{k},z_{k+1},\dots, z_{n})dz_{k},\quad k=1,\dots, n-1,
\end{align*}
where $\Gamma_{z}^{(k)}$ is a counterclockwise integral contour enclosing the points $z\pm k$ and $z\pm 1$,
but no other poles of the integrand.
Finally, we set
\begin{align*}
	I^{(n)}_{\lambda} = \int_{\Gamma}I_{\lambda}^{(n-1)}(z)dz,
\end{align*}
where $\Gamma$ is a counterclockwise integral contour enclosing all poles of $I_{\lambda}^{(n-1)}$ on $\bC$.
Then, we have
\begin{equation*}
	\Sigma_{(n)}(\lambda) = -\frac{1}{k}I^{(n)}_{\lambda}.
\end{equation*}
\end{thm}

\begin{figure}
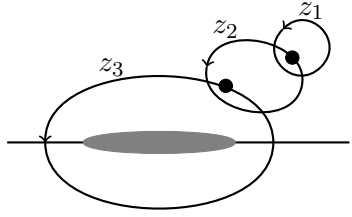

\ctikzfig{satellite}
\caption{Satellite integral for $n=3$. Suppose that the poles from $H_{\lambda}$ are in the grey region. Each black dot actually consists of a few poles.}
\label{fig:satellite}
\end{figure}

Before going for the proof, let us evoke a visual picture for $I^{(n)}_{\lambda}$ as a {\it satellite integral} of the function $F^{(n)}_{\lambda}$ (Figure~\ref{fig:satellite}).

\begin{proof}
Due to the Frobenius character formula (\ref{eq:FrobeniusCharFormula}), we need to show
\begin{align}
\label{eq:intISatelite}
	\int_{\Gamma} I_{\lambda}^{(n-1)}(z) dz = \int_{\Gamma} H_{\lambda}(z) H_{\lambda}(z-1)\cdots H_{\lambda}(z-n+1) dz.
\end{align}
For that, we first show
\begin{align}
\label{eq:formulaForI}
	 &\, I_{\lambda}^{(k)}(z_{k+1},\dots, z_{n}) \\
	=&\, \frac{1}{2}\Bigg( \frac{H_{\lambda}(z_{k+1})H_{\lambda}(z_{k+1}+1)\cdots H_{\lambda}(z_{k+1}+k)}{z_{k+1}-z_{k+2}}
	\prod_{j=k+2}^{n}\frac{(z_{k+1}-z_{j})(z_{k+1}-z_{j}+k)}{(z_{k+1}-z_{j}-1)(z_{k+1}-z_{j}+k+1)} \notag \\
	&\, +\frac{H_{\lambda}(z_{k+1})H_{\lambda}(z_{k+1}-1)\cdots H_{\lambda}(z_{k+1}-k)}{z_{k+1}-z_{k+2}}
	\prod_{j=k+2}^{n}\frac{(z_{k+1}-z_{j})(z_{k+1}-z_{j}-k)}{(z_{k+1}-z_{j}+1)(z_{k+1}-z_{j}-k-1)}\Bigg) \notag \\
	&\, \cdot F_{\lambda}^{(n-k-1)}(z_{k+2},\dots,z_{n})\notag
\end{align}
for $k=0,\dots,n-2$ by induction on $k$.
When $k=0$, (\ref{eq:formulaForI}) is clear from
\begin{equation}
\label{eq:functionFReduction}
	F_{\lambda}^{(n)}(z_{1},\dots,z_{n}) = \frac{H_{\lambda}(z_{1})}{z_{1}-z_{2}}\prod_{i=2}^{n}\frac{(z_{1}-z_{i})^{2}}{(z_{1}-z_{i}-1)(z_{1}-z_{i}+1)}\, \cdot F_{\lambda}^{(n-1)}(z_{2},\dots, z_{n}).
\end{equation}
Let us assume (\ref{eq:formulaForI}) for some $k$ and prove the case of $k+1$.
We are integrating out $z_{k+1}$ along $\Gamma^{(k+1)}_{z_{k+2}}$, we are also allowed to translate the integral variable.
In the second term of the sum in the right-hand side, we change the integral variable $z_{k+1}$ to $z_{k+1}+k$ and get
\begin{align*}
	&\, I_{\lambda}^{(k+1)}(z_{k+2},\dots, z_{n}) \\
	=&\, \frac{1}{2}\int_{\Gamma^{(k+1)}_{z_{k+2}}}\Bigg(\frac{1}{z_{k+1}-z_{k+2}}+\frac{1}{z_{k+1}-z_{k+2}+k}\Bigg)
	H_{\lambda}(z_{k+1})\cdots H_{\lambda}(z_{k+1}+k) \\
	&\qquad \cdot \prod_{j=k+2}^{n}\frac{(z_{k+1}-z_{j})(z_{k+1}-z_{j}+k)}{(z_{k+1}-z_{j}-1)(z_{k+1}-z_{j}+k+1)}
	F_{\lambda}^{(n-k-1)}(z_{k+2},\dots,z_{n})dz_{k+1} \\
\end{align*}
In the integral, we can perform the following rearrangement:
\begin{align*}
	&\, \Bigg(\frac{1}{z_{k+1}-z_{k+2}}+\frac{1}{z_{k+1}-z_{k+2}+k}\Bigg)\frac{(z_{k+1}-z_{k+2})(z_{k+1}-z_{k+2}+k)}{(z_{k+1}-z_{k+2}-1)(z_{k+1}-z_{k+2}+k+1)} \\
	=&\, \frac{2(z_{k+1}-z_{k+2})+k}{(z_{k+1}-z_{k+2}-1)(z_{k+1}-z_{k+2}+k+1)} \\
	=&\frac{1}{z_{k+1}-z_{k+2}-1}+\frac{1}{z_{k+1}-z_{k+2}+k+1}.
\end{align*}
Thus we obtain
\begin{align*}
	&\, I_{\lambda}^{(k+1)}(z_{k+2},\dots, z_{n}) \\
	=& \frac{1}{2}\int_{\Gamma^{(k+1)}_{z_{k+2}}}\Bigg(\frac{1}{z_{k+1}-z_{k+2}-1}+\frac{1}{z_{k+1}-z_{k+2}+k+1}\Bigg)
	H_{\lambda}(z_{k+1})\cdots H_{\lambda}(z_{k+1}+k) \\
	&\qquad \cdot \prod_{j=k+3}^{n}\frac{(z_{k+1}-z_{j})(z_{k+1}-z_{j}+k)}{(z_{k+1}-z_{j}-1)(z_{k+1}-z_{j}+k+1)}
	F_{\lambda}^{(n-k-1)}(z_{k+2},\dots,z_{n})dz_{k+1} \\
	=&\frac{1}{2}\Bigg( H_{\lambda}(z_{k+2}+1)\cdots H_{\lambda}(z_{k+2}+k+1)
	 \prod_{j=k+3}^{n}\frac{(z_{k+2}-z_{j}+1)(z_{k+2}-z_{j}+k+1)}{(z_{k+2}-z_{j})(z_{k+2}-z_{j}+k+2)} \\
	&\, + H_{\lambda}(z_{k+2}-1)\cdots H_{\lambda}(z_{k+2}-k-1)
	 \prod_{j=k+3}^{n}\frac{(z_{k+2}-z_{j}-k-1)(z_{k+2}-z_{j}-1)}{(z_{k+2}-z_{j}-k-2)(z_{k+2}-z_{j})}\Bigg) \\
	&\quad \cdot F_{\lambda}^{(n-k-1)}(z_{k+2},\dots,z_{n}) \\
	=&\, \frac{1}{2}\Bigg( \frac{H_{\lambda}(z_{k+2})\cdots H_{\lambda}(z_{k+2}+k+1)}{z_{k+2}-z_{k+3}}
	\prod_{j=k+3}^{n}\frac{(z_{k+2}-z_{j})(z_{k+2}-z_{j}+k+1)}{(z_{k+2}-z_{j}-1)(z_{k+2}-z_{j}+k+2)} \\
	&\, + \frac{H_{\lambda}(z_{k+2})\cdots H_{\lambda}(z_{k+2}-k-1)}{z_{k+2}-z_{k+3}}
	 \prod_{j=k+3}^{n}\frac{(z_{k+2}-z_{j})(z_{k+2}-z_{j}-k-1)}{(z_{k+2}-z_{j}+1)(z_{k+2}-z_{j}-k-2)}\Bigg) \\
	&\quad \cdot F_{\lambda}^{(n-k-2)}(z_{k+3},\dots,z_{n}).
\end{align*}
To get the last line, we used (\ref{eq:functionFReduction}).
Thus, (\ref{eq:formulaForI}) has been proven. In particular, we have
\begin{align*}
	&\, I^{(n-2)}_{\lambda}(z_{n-1},z_{n}) \\
	=&\, \frac{1}{2}\Bigg(H_{\lambda}(z_{n-1})\cdots H_{\lambda}(z_{n-1}+n-2)\frac{z_{n-1}-z_{n}+n-2}{(z_{n-1}-z_{n}-1)(z_{n-1}-z_{n}+n-1)} \\
	&\, + H_{\lambda}(z_{n-1})\cdots H_{\lambda}(z_{n-1}-n+2)\frac{z_{n-1}-z_{n}-n+2}{(z_{n-1}-z_{n}+1)(z_{n-1}-z_{n}-n+1)}\Bigg)H_{\lambda}(z_{n}).
\end{align*}
We integrate this once more over $z_{n-1}$ along $\Gamma^{(n-1)}_{z_{n}}$.
The actual computation is similar to the previous one; after change of the integration variable, we get
\begin{align*}
	I_{\lambda}^{(n-1)}(z) =& \frac{1}{2}\Big( H_{\lambda}(z)H_{\lambda}(z+1)\cdots H_{\lambda}(z+n-1)
	+H_{\lambda}(z)H_{\lambda}(z-1)\cdots H_{\lambda}(z-n+1)\Big).
\end{align*}
Finally, integrating over $z$ along $\Gamma$, we get the desired (\ref{eq:intISatelite}).
\end{proof}

\subsection{Diagram elements as integral}
\begin{figure}
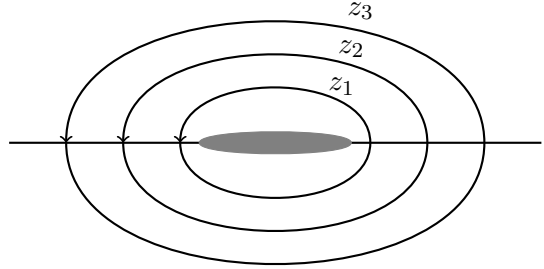

\ctikzfig{radialIntegral}
\caption{Radial integral for $\cI^{(3)}_{\lambda}[\id]$. The integrand is the same as Figure~\ref{fig:satellite}.}
\label{fig:radialIntegral}
\end{figure}

Next, we see that the right-hand side of (\ref{eq:characterDiagramEval}) can be expressed as an integral of
the same function as above, but along different contours.
Let $\lambda\in \bY$ and $n\in\bZ_{>0}$.
For each $\sigma\in S_{n}$, we write $\Gamma^{(n)}_{\lambda} [\sigma]$ for integral contours for $n$ variables $z_{1},\dots,z_{n}$
such that
\begin{itemize}
\item 	$|z_{\sigma (i+1)}|>|z_{\sigma_{i}}\pm 1|$ for all $i=1,\dots, n-1$,
\item 	each variable goes counterclockwise,
\item 	$z_{\sigma (1)}$ encloses all poles of $H_{\lambda}$,
\end{itemize}
and set
\begin{align*}
	\cI^{(n)}_{\lambda} [\sigma] = \int_{\Gamma^{(n)}_{\lambda} [\sigma]}F_{\lambda}^{(n)}(z_{1},\dots,z_{n})dz_{1}\cdots dz_{n},
\end{align*}
where the integrand is the same one given in (\ref{eq:defFunctionF}).
This time, $\cI^{(n)}_{\lambda}$ was given by a radial integral (Figure~\ref{fig:radialIntegral}).

\begin{thm}
\label{thm:DiagramCharacterIntegral}
For $\lambda\in \bY$, we have
\begin{align*}
	\cI^{(n)}_{\lambda} [\id]=(-1)^{n}\, \input{characterOneCycleEval.tikz}\, .
\end{align*}
\end{thm}
\begin{proof}
Let us first see what happens when we integrate out $z_{1}$ from $F_{\lambda}^{(n)}(z_{1},\dots,z_{n})$ along a contour that encloses all poles of $H_{\lambda}(z_{1})$: for each pole of $H_{\lambda}(z_{1})$, we find a residue
\begin{align*}
	-|\lambda|p^{\downarrow}(\lambda,\lambda^{(1)}) = -|\lambda|\frac{\dim\lambda^{(1)}}{\dim \lambda}
\end{align*}
with some $\lambda^{(1)}$ such that $\lambda^{(1)}\nearrow\lambda$ (see (\ref{eq:cotrans_Cauchy_contents}) and \ref{eq:cotrans_dim_ratio})).
At the same time, for each $i=2,\dots, n$, we have the formula
\begin{align*}
	\frac{(c(\lambda/\lambda^{(1)})-z_{i})^{2}}{(c(\lambda/\lambda^{(1)})-z_{i}-1)(c(\lambda/\lambda^{(1)})-z_{i}+1)}H_{\lambda}(z_{i})=H_{\lambda^{(1)}}(z_{i})
\end{align*}
from (\ref{eq:Cauchy_transform_adding_box}).
Thus, we observe
\begin{align*}
	\int_{\Gamma_{\lambda}^{(1)}[\id]}F_{\lambda}^{(n)}(z_{1},\dots,z_{n})dz_{1}
	= -\frac{|\lambda|}{\dim \lambda}\sum_{\lambda^{(1)};\lambda^{(1)}\nearrow \lambda} \dim \lambda^{(1)}\frac{1}{c(\lambda/\lambda^{(1)})-z_{2}}F^{(n-1)}_{\lambda^{(1)}}(z_{2},\dots,z_{n}).
\end{align*}

Repeating this calculation, we obtain the identities
\begin{align*}
	&\, \int_{\Gamma^{(k)}_{\lambda}[\id]}F_{\lambda}^{(n)}(z_{1},\dots,z_{n})dz_{1}\cdots dz_{k} \\
	=&\, (-1)^{k}\frac{(|\lambda|\downharpoonright k)}{\dim \lambda}\sum_{\substack{\lambda^{(1)},\dots,\lambda^{(k)}\\ \lambda^{(k)}\nearrow \cdots \nearrow \lambda^{(1)}\nearrow \lambda}}\dim \lambda^{(k)}
	\prod_{i=1}^{k-1}\frac{1}{c(\lambda^{(i-1)}/\lambda^{(i)})-c(\lambda^{(i)}/\lambda^{(i+1)})}\\
	&\qquad \cdot \frac{1}{c(\lambda^{(k-1)}/\lambda^{(k)})-z_{k+1}}\, F_{\lambda^{(k)}}(z_{k+1},\dots,z_{n})
\end{align*}
for $k=1,\dots, n$.
In particular, we can see that the case of $k=n$ gives the desired result comparing it with (\ref{eq:characterDiagramEval}).
\end{proof}

\subsection{Change of integral contours}

Combining Theorems~\ref{thm:FrobeniusSatelite} and \ref{thm:DiagramCharacterIntegral}, we can state Theorem~\ref{thm:characterDiagram} in the following way.

\begin{thm}
For $\lambda\in \bY$ and $n\in \bZ_{>0}$, we have
\begin{equation*}
    (-1)^{n}\cI_{\lambda}^{(n)}[\id] = -\frac{1}{n}I^{(n)}_{\lambda}.
\end{equation*}
\end{thm}

We have proved this by representation theory, but it could be possible to prove this directly since both sides are integrals of the same function $F_{\lambda}^{(n)}$ along different contours.
We have not managed this, but can only display computations for small $n$.
If it succeeds for all $n$, we can also say that we get a direct proof of the Frobenius formula that does not bypass symmetric functions.

The following properties will be useful for general $n$.
\begin{lem}
\label{lem:functionFcyclic}
Let $\lambda\in \bY$ and $n\in\bZ_{>0}$.
Then, we have
\begin{align*}
	\sum_{p=1}^{n}F_{\lambda}^{(n)}(z_{p},z_{p+1} \dots,z_{n},z_{1},z_{2}\dots ,z_{p-1}) = 0.
\end{align*}
\end{lem}
\begin{proof}
Notice that
\begin{equation}
\label{eq:functionFSymmetric}
	\prod_{i=1}^{n-1}(z_{i}-z_{i+1}) F_{\lambda}^{(n)}(z_{1},\dots, z_{n})
\end{equation}
is a symmetric function of $z_{1},\dots, z_{n}$.
In particular, we have
\begin{align*}
	(z_{n}-z_{1}) F_{\lambda}^{(n)}(z_{p},\dots, z_{n},z_{1},\dots, z_{p-1}) = (z_{p-1}-z_{p})F_{\lambda}^{(n)}(z_{1},\dots, z_{n})
\end{align*}
for $p=2,\dots, n$.
Thus, we get
\begin{align*}
	\sum_{p=1}^{n}F_{\lambda}^{(n)}(z_{p},\dots, z_{n},z_{1},\dots, z_{p-1}) = \Big(1 + \sum_{p=2}^{n-1}\frac{z_{p-1}-z_{p}}{z_{n}-z_{1}}\Big) F_{\lambda}^{(n)}(z_{1},\dots, z_{n}) = 0
\end{align*}
as was desired.
\end{proof}

\begin{lem}
\label{lem:functionFinversion}
Let $\lambda\in \bY$ and $n\in\bZ_{>0}$.
Then, we have
\begin{align*}
	F_{\lambda}^{(n)}(z_{n},z_{n-1},\dots, z_{1}) = (-1)^{n-1}F_{\lambda}^{(n)}(z_{1},z_{2}, \dots, z_{n}).
\end{align*}
\end{lem}
\begin{proof}
The proof again uses the fact that (\ref{eq:functionFSymmetric}) is symmetric.
Indeed, the product $\prod_{i=1}^{n-1}(z_{i}-z_{i+1})$ yields the sign $(-1)^{n-1}$ after inversion.
\end{proof}

\subsubsection{The $n=2$ case}
\begin{figure}
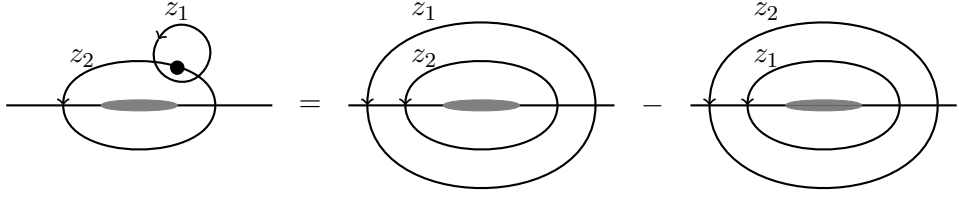

\ctikzfig{changeOfContours}
\caption{Change of integral contours.}
\label{fig:changeOfContours}
\end{figure}

By deforming the integral contours in $I^{(2)}_{\lambda}$ (see Figure~\ref{fig:changeOfContours}), we get
\begin{align*}
	I^{(2)}_{\lambda} = \cI^{(2)}_{\lambda}[(2\,1)] - \cI^{(2)}_{\lambda}[\id].
\end{align*}
As a notation, we write $(\sigma (1)\,\dots\,\sigma (n))$ to express a permutation $\sigma$.
In this case of $n=2$, Lemmas~\ref{lem:functionFcyclic} and \ref{lem:functionFinversion} give the same result
\begin{align*}
	\cI^{(2)}_{\lambda}[(2\, 1)] = -\cI^{(2)}_{\lambda}[\id]
\end{align*}
after change of integration variables.
Therefore, the identity
\begin{equation*}
	\cI^{(2)}_{\lambda}[\id] = -\frac{1}{2}I^{(2)}_{\lambda}
\end{equation*}
holds.

\subsubsection{The $n=3$ case}
When we similarly deform the integral contours, we get
\begin{align*}
	I^{(3)}_{\lambda}
	&= \cI_{\lambda}^{(3)}[(3\,2\,1)] +  \cI_{\lambda}^{(3)}[(1\,2\,3)] -  \cI_{\lambda}^{(3)}[(3\,1\,2)] -  \cI_{\lambda}^{(3)}[(2\,1\,3)] \\
	&= 2 \cI_{\lambda}^{(3)}[\id] - (\cI_{\lambda}^{(3)}[(3\,1\,2)] +  \cI_{\lambda}^{(3)}[(2\,1\,3)]).
\end{align*}
Here, we used Lemma~\ref{lem:functionFinversion} to get the second line.

\begin{lem}
We have $\cI_{\lambda}^{(3)}[(2\,1\,3)] = \cI_{\lambda}^{(3)}[(2\,3\,1)]$.
\end{lem}
\begin{proof}
The difference $\cI_{\lambda}^{(3)}[(2\,1\,3)] - \cI_{\lambda}^{(3)}[(2\,3\,1)]$ is
\begin{align*}
	\cI_{\lambda}^{(3)}[(2\,1\,3)] - \cI_{\lambda}^{(3)}[(2\,3\,1)] = \int \Big( \int \Big(\int_{C}F^{(3)}(z_{1},z_{2},z_{3})dz_{1}\Big) dz_{3}\Big)dz_{2},
\end{align*}
where $C$ encloses $z_{3}\pm 1$ but no other poles, $z_{3}$ encloses all poles of $H_{\lambda}$ and $z_{2}\pm 1$, and $z_{2}$ encloses all poles of $H_{\lambda}$.
Let us set
\begin{align*}
    \wtilde{F}_{\lambda}^{(3)}(z_{1},z_{2},z_{3}) = \frac{(z_{1}-z_{3}-1)(z_{1}-z_{3}+1)}{(z_{1}-z_{3})^{2}}F_{\lambda}^{(3)}(z_{1},z_{2},z_{3}).
\end{align*}
When we integrate $F_{\lambda}^{(3)}$ in $z_{1}$ along a contour $C$, we collects the contribution from the poles in
\begin{align*}
    \frac{(z_{1}-z_{3})^{2}}{(z_{1}-z_{3}-1)(z_{1}-z_{3}+1)} =\frac{1}{2}\Big(\frac{z_{1}-z_{3}}{z_{1}-z_{3}-1}+\frac{z_{1}-z_{3}}{z_{1}-z_{3}+1}\Big).
\end{align*}
Thus, we get
\begin{align*}
    \int_{C}F_{\lambda}^{(3)}(z_{1},z_{2} ,z_{3})dz_{1}
    =\frac{1}{2}\Big(\wtilde{F}_{\lambda}^{(3)}(z_{3}+1,z_{2}, z_{3}) - \wtilde{F}_{\lambda}^{(3)}(z_{3}-1,z_{2}, z_{3})\Big).
\end{align*}
It follows from the symmetry of (\ref{eq:functionFSymmetric}) that $\wtilde{F}^{(3)}(z_{1},z_{2},z_{3})$ is symmetric under permuting $z_{1}$ and $z_{3}$. 
Hence, we have
\begin{align*}
    \int_{C}F_{\lambda}^{(3)}(z_{1},z_{2} ,z_{3})dz_{1}
    =\frac{1}{2}\Big(\wtilde{F}_{\lambda}^{(3)}(z_{3}+1,z_{2}, z_{3}) - \wtilde{F}_{\lambda}^{(3)}(z_{3},z_{2}, z_{3}-1)\Big),
\end{align*}
but as long as we integrate this over $z_{3}$, the two terms above cancel each other:
\begin{align*}
	\int \Big(\int_{C}F_{\lambda}^{(3)}(z_{1},z_{2} ,z_{3})dz_{1}\Big) dz_{3} = 0.
\end{align*}
Therefore, we have $\cI_{\lambda}^{(3)}[(2\,1\,3)] = \cI_{\lambda}^{(3)}[(2\,3\,1)]$.
\end{proof}

The above lemma allows us to get
\begin{align*}
	I^{(3)}_{\lambda}
	&= 2 \cI_{\lambda}^{(3)}[\id] - (\cI_{\lambda}^{(3)}[(3\,1\,2)] +  \cI_{\lambda}^{(3)}[(2\,3\,1)]) = 3 \cI_{\lambda}^{(3)}[\id],
\end{align*}
where we also used Lemma~\ref{lem:functionFcyclic}.
Therefore, we have
\begin{equation*}
	- \cI_{\lambda}^{(3)}[\id] = -\frac{1}{3}I^{(3)}_{\lambda}.
\end{equation*}

\bibliographystyle{alpha}
\bibliography{heis_cat}

\end{document}